\documentstyle[12pt,psfig]{article}

\def\hatM{{\widehat{M}}}
\def\hatv{{\widehat{v}}}
\def\hattt{{\widehat{\calt}}}
\def\caltbar{\overline{\calt}}
\def\hatpp{{\widehat{\pp}}}
\def\interior{{\rm Int}}
\def\index{{\rm ind}}
\def\vbar{{\overline v}}
\def\Wbar{{\overline W}}
\def\Kbar{{\overline K}}

\def\hbar{{\overline h}}
\def\GL{{\rm GL}}

\def\Ccell{{\rm C}^{\rm cell}}
\def\Cphi{{\rm C}^\varphi}

\def\basis{{\hbox{\Got b}}}
\def\hbasis{{\hbox{\Got h}}}
\def\gbasis{{\hbox{\Got g}}}

\def\ristr#1{\big|_{#1}}

\def\stwotriv{S^2_{\rm triv}}

\reversemarginpar

\makeatletter
\def\@begintheorem#1#2{\it \trivlist \item[\hskip \labelsep{\bf #1\ #2.}]}
\newtheorem{teo}{Theorem}[section]
\newtheorem{rem}[teo]{Remark}
\newtheorem{lem}[teo]{Lemma}

\newtheorem{prop}[teo]{Proposition}

\def\finedim#1{{\hfill\hbox{\enspace\fbox{\ref{#1}}}}\vspace{5pt}}
\def\dim#1{\vspace{1pt}\noindent{\it Proof of} {\hspace{2pt}}\ref{#1}.}
\def\compo{\,{\scriptstyle\circ}\,}
\def\cont{{\rm C}}

\font\scpicc=cmcsc10
\font\sc=cmcsc10 scaled 1200

\newfont{\Bbb}{msbm10 scaled 1200}
\def\mr{{\hbox{\Bbb R}}}

\def\mz{{\hbox{\Bbb Z}}}

\newfont{\mycal}{eusm10 scaled 1200}

\newfont{\Got}{eufm10 scaled 1200}

\def\cc{{\cal C}}

\def\uu{{\cal U}}
\def\ss{{\cal S}}

\def\calt{{\cal T}}

\def\pp{{\cal P}}

\def\dd{{\cal D}}

\def\eul{{\rm Eul}}
\def\euls{{\rm Eul}^{\rm s}}
\def\eulc{{\rm Eul}^{\rm c}}
\def\Thetas{\Theta^{\rm s}}
\def\Thetac{\Theta^{\rm c}}
\def\alphas{\alpha^{\rm s}}
\def\alphac{\alpha^{\rm c}}

\title{Reidemeister Torsion of 3-Dimensional\\ Euler Structures
with Simple Boundary Tangency\\ and Legendrian Knots}

\author{Riccardo Benedetti\qquad Carlo Petronio\thanks{The second named
author gratefully acknowledges financial support by GNSAGA-CNR}}

\begin{document}

\maketitle

\noindent{\small{\scpicc Abstract}. We generalize Turaev's definition of torsion invariants of pairs
$(M,\xi)$, where $M$ is a 3-dimensional manifold and $\xi$ is an Euler structure on $M$ (a non-singular
vector field up to homotopy relative to $\partial M$  and local modifications in $\interior(M)$). Namely, we
allow $M$ to have arbitrary boundary and $\xi$ to have simple (convex and/or concave) tangency circles to
the boundary.  We prove that Turaev's $H_1(M)$-equivariance formula holds also in  our generalized context.
Our torsions apply in particular to (the exterior of) Legendrian links (in particular, knots) in contact
3-manifolds, and we  prove that they can distinguish knots which are isotopic as framed knots but not as
Legendrian knots. Using the combinatorial encoding of vector fields based on branched standard spines we
show how to explicitly invert Turaev's reconstruction map  from combinatorial to smooth Euler structures,
thus making the computation of torsions a more effective one. As an example we work out a specific
computation.}

\vspace{.5cm}

\noindent{\small{\scpicc Mathematics Subject Classification (1991)}: 57N10
(primary), 57Q10, 57R25 (secondary).}

\vspace{.5cm}

\section*{Introduction}
Reidemeister torsion is a classical yet very vital topic in
3-dimensional topology, and it was recently used in a variety of important
developments. To mention a few, torsion is a fundamental ingredient of the
Casson-Walker-Lescop invariants (see {\em e.g.}~\cite{lescop}), and more
generally of the perturbative approach to quantum invariants (see {\em
e.g.}~\cite{lemuoh}). Relations have been pointed out between torsion and
hyperbolic geometry~\cite{porti}. Turaev's torsion of Euler 
structures~\cite{turaev:Euler} has recently been 
recognized by Turaev himself (\cite{turaev:spinc},~\cite{turaev:nuovo}) to have deep
connections with the Seiberg-Witten
invariants of ${\rm Spin}^{\rm c}$-structures on 3-manifolds, after 
the proof of Meng and Taubes~\cite{meng:taub} that a suitable combination of these
invariants can be identified with the classical Milnor torsion.

Turaev's theory~\cite{turaev:Euler} actually exists in all dimensions. We
quickly review it before proceeding. A {\em smooth Euler structure} $\xi$ on a
compact oriented\footnote{Orientability is  not strictly necessary, but we find
it convenient to assume it.} manifold $M$, possibly with $\partial
M=\emptyset$, is a non-singular vector field on $M$ viewed up to local
modifications  in $\interior(M)$ and  homotopy relative to $\partial M$. Turaev
allows only ``monochromatic'' boundary components, {\it i.e.}~black ones (on
which the field points outwards) and white ones (on which it points inwards).
This implies the constraint that $\chi(M,W)=0$, where $W$ is the white portion
of  $\partial M$, but in~\cite{turaev:spinc} and~\cite{turaev:nuovo} Turaev
only focuses on the more specialized case where $M$ is 3-dimensional and closed
or bounded by tori. In all dimensions, the set $\euls(M,W)$ of smooth Euler
structures compatible with $(M,W)$ is an affine space over  $H_1(M;\mz)$. The
two main ingredients of Turaev's theory are as follows. First, he defines a
certain set of $1$-chains, called the space $\eulc(M,W)$ of {\em combinatorial}
Euler structures compatible with $(M,W)$, he shows that this is again affine
over $H_1(M;\mz)$,  and he describes an $H_1(M;\mz)$-equivariant bijection
$\Psi:\eulc(M,W)\to\euls(M,W)$ called the {\em reconstruction map}. Second, for
$\xi\in\eulc(M,W)$ and for  any representation $\varphi$ of $\pi_1(M)$ into the
units of a suitable ring $\Lambda$ he defines a torsion  invariant
$\tau^\varphi(M,\xi)$, or more generally $\tau^\varphi(M,\xi,\hbasis)$, with
values in $K_1(\Lambda)/(\pm1)$. This invariant is by definition a lifting of
the classical Reidemeister torsion (see~\cite{milnor}) $\tau^\varphi(M)\in
K_1(\Lambda)/(\pm\varphi(\pi_1(M)))$, and it satisfies the
$H_1(M;\mz)$-equivariance formula \begin{equation} \label{tau:equivariance}
\tau^\varphi(M,\xi',\hbasis)=\tau^\varphi(M,\xi,\hbasis)\cdot \varphi(\xi'-\xi)
\end{equation} where $\xi'-\xi\in H_1(M;\mz)$. For $\xi\in\euls(M,W)$ one
defines $\tau^\varphi(M,\xi)$ as $\tau^\varphi(M,\Psi^{-1}(\xi))$, and the 
$H_1(M;\mz)$-equivariance of the reconstruction map $\Psi$ implies that 
formula~(\ref{tau:equivariance}) holds also for smooth structures. We emphasize
that the definition of $\Psi$ is based on an explicit geometric construction,
but its bijectivity is only established through $H_1(M;\mz)$-equivariance. This
makes the definition of torsion for smooth structures somewhat implicit.

In the present paper, and in other papers in preparation, we are concerned with
generalizations and improvements of Turaev's theory. Here we consider
3-manifolds. This work had two main initial aims. Our first aim was to find a
geometric description of the map $\Psi^{-1}$, and hence to turn the computation
of Turaev's torsion into a more effective procedure, using our combinatorial
encoding~\cite{lnm} of non-singular vector fields up to homotopy (also
called ``combings'') in terms of branched standard spines. Our second aim was
to define torsion invariants of (pseudo-)Legendrian links, {\em i.e.}~links
tangent to a given plane field, viewed up to tangency-preserving
isotopy (when the plane field is a contact
structure one gets the familiar notion of Legendrian link and Legendrian
isotopy). A specific motivation to look for invariants of pseudo-Legendrian
links comes  from the remarkable relation recently discovered by Fintushel and
Stern~\cite{fist} between  the Alexander polynomial ({\em i.e.}~Milnor torsion)
of a knot $K\subset S^3$ and  (a suitable combination of) the Seiberg-Witten
invariants of the  ``surgered'' $4$-manifold $X_K$ obtained using $K$ (and a
suitable base $4$-manifold $X$). Both our initial aims lead us to consider
Euler structures on $3$-manifolds $M$ (without restrictions on $\partial M$)
allowing simple tangency circles to $\partial M$ of  {\em concave} type (see
Fig.~\ref{conc:conv:tang} below). On the other  hand it turns out that, to
define torsions, the natural objects to deal with  are Euler structures with
{\em convex} tangency  circles. It is a fortunate fact, peculiar of dimension
3, that there is a canonical way to associate  a convex field to any {\em
simple} ({\em i.e.}~mixed concave and convex) one. This allows to define
torsion for all smooth simple Euler structures, and eventually to achieve both
the objectives we had in mind.

Let us now summarize the contents of this paper. The foundational part of
our work consists in extending to the
context of Euler structures with simple tangency the notions of combinatorial
structure $\eulc$ and reconstruction map $\Psi$. 
This part follows the same scheme as~\cite{turaev:Euler}
and relies on technical results of Turaev. Our main contribution here
is the proof that the natural transformations of a concave structure
into a convex one, viewed at the smooth level and at the combinatorial level,
actually correspond to each other under the reconstruction map
(Theorem~\ref{diagram:commutes}).
After setting the foundations, we prove the following main results
(stated informally here:
see Sections~\ref{statements} 
and~\ref{link:section} for precise definitions and statements.)

\begin{teo}\label{informal:good:for:links}
There exist pairs $(M,\eta)$, where $\eta$ is an oriented plane distribution on $M$,
and knots $K_0$, $K_1$ tangent to $\eta$ (whence framed), such that:
\begin{enumerate}
\item $K_0$ and $K_1$ are isotopic as framed knots;
\item A torsion invariant shows that $K_0$ and $K_1$ are not isotopic to each
other through knots tangent to $\eta$.
\end{enumerate}
\noindent Moreover, $\eta$ can be chosen to be a contact structure.
\end{teo}

\begin{teo}\label{informal:algorithmic}
Let $\xi$ be an Euler structure with concave
tangency circles. If $P$ is a branched standard spine 
which represents $\xi$, then $P$ allows to explicitly 
find a representative of $\Psi^{-1}(\xi)\in\eulc$, and hence to compute
the torsion of $\xi$ in terms of the finite combinatorial data which encode $P$.
\end{teo}

Concerning Theorem~\ref{informal:good:for:links}, 
the relation with vector fields is of course given by taking the orthogonal to $\eta$. 
Note that our torsions are defined not only for knots but also for links, including
homologically non-trivial Legendrian links, for which the usual invariants, such
as the rotation number (Maslov index), are not defined. However 
we will show that, in a suitable sense, torsion is a generalization of the rotation number,
when the latter is defined. Moreover we will prove
that torsions are sensitive to an analogue of the winding number.

Another relevant point concerning
Theorem~\ref{informal:good:for:links} is that in general torsion does not
provide a single-valued invariant for pseudo-Legendrian knots, because the
action of a certain mapping class group must be taken into account. However we
will show that for many knots this action is actually trivial, so torsion is
indeed single-valued. This is the case for instance for all knots contained in
homology spheres. So our torsion invariants include a (non-trivial) 
refinement of the Alexander polynomial to Legendrian knots in $S^3$.

This paper is organized as follows. In Section~\ref{statements} we provide the
formal definitions of smooth and combinatorial Euler structures, we informally
introduce torsion and we state the fundamental results which imply that 
torsion is well-defined and 
$H_1$-equivariant. In Section~\ref{link:section} we specialize to knot exteriors and
we prove that torsion is a non-trivial invariant of Legendrian knots, as announced above.
In Section~\ref{torsion:def:section} we address the
main technical points of the definition of torsion. In Section~\ref{spines:section} we
show how branched standard spines can be used for computing torsion, and in Section~\ref{exa:section} we actually carry out a computation. In Sections~\ref{statements}
to~\ref{spines:section} proofs which are long and require the introduction of 
ideas and techniques not used elsewhere are omitted. Section~\ref{proofs} contains all
these proofs.

We conclude this introduction by announcing related results which we have
recently obtained and currently writing down. In~\cite{second:paper} we extend
to the case with boundary our combinatorial presentation~\cite{lnm} of combed
manifolds in terms of branched spines (this extension is actually mentioned
also in the present paper ---see  Section~\ref{spines:section}). Using this
presentation we develop an approach to torsion entirely based on combinatorial
techniques (branched spines),  leading to a generalization of Turaev's theory
slightly different from the present paper's. In~\cite{third:paper} we
generalize the theory of Euler structures and (with some restrictions) of
torsions to all dimensions and allowing any generic (Whitney-Morin-type)
tangency to the boundary.  Noting that this situation arises when one cuts a
manifold along a hypersurface in general position with respect to a given
non-singular vector field, one is naturally lead to the question of how Euler
structures and torsion behave under glueing. As a motivation, note that a
similar question is involved in the summation formulae for the Casson invariant
(see~\cite{lescop:bis}), and is faced also 
in~\cite{fist},~\cite{meng:taub} and~\cite{turaev:nuovo}. We believe
that this question deserves further investigation.

\section{Main definitions and statements}\label{statements}

In this section we define Euler structures and their torsion.
Fix once and for ever a compact oriented 3-manifold $M$, possibly with 
$\partial M=\emptyset$.
Using the {\it Hauptvermutung}, we will always freely intermingle the
differentiable, piecewise linear and topological viewpoints. Homeomorphisms
will always respect orientations. All vector fields mentioned in this paper will be 
non-singular, and they will be termed just {\em fields} for the sake of brevity.

\begin{figure}
\centerline{\psfig{file=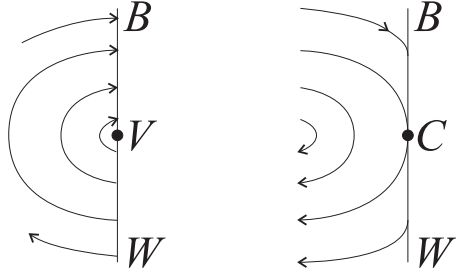,width=6cm}}
\caption{\label{conc:conv:tang} Convex (left) and concave (right) tangency to the boundary.}
\end{figure}

\paragraph{Smooth and combinatorial Euler structures.}
We will call {\em boundary pattern} on $M$ a partition $\pp=(W,B,V,C)$ of
$\partial M$ where $V$ and $C$ are finite unions of 
disjoint circles, and $\partial W=\partial B=V\cup C$. In particular, $W$ and $B$ are
interiors of compact surfaces embedded in $\partial M$. Even if $\pp$ can actually 
be determined by less data, {e.g.} the pair $(W,V)$,
we will find it convenient to refer to $\pp$ as a quadruple.
Points of $W$, $B$, $V$ and $C$ will be called {\em white}, {\em black}, 
{\em convex} and {\em concave} respectively.
We define the set
of {\em smooth Euler structures} on $M$ compatible with $\pp$, denoted by
$\euls(M,\pp)$, as the set of equivalence classes of
fields on $M$ which point inside on $W$, point
outside on $B$ and have simple tangency to $\partial M$ 
of {\em convex} type along $V$ and {\em concave} type 
along $C$,
as shown in a cross-section in Fig.~\ref{conc:conv:tang}.
Two such fields are equivalent if they are obtained from each other by 
homotopy through fields of the same type and modifications supported into
interior balls. The following
variation on the Poincar\'e-Hopf formula is established in Section~\ref{proofs}:

\begin{prop}\label{p:h:formula}
$\euls(M,\pp)$ is non-empty if and only if $\chi(\Wbar)=\chi(M)$. 
\end{prop}

\noindent We remark here that $\chi(\Wbar)=\chi(W)$, $\chi(\overline{B})=\chi(B)$, $\chi(V)=\chi(C)=0$ and
$\chi(W)+\chi(B)=\chi(\partial M)=2\chi(M)$, so there are various ways to
rewrite the relation $\chi(\Wbar)=\chi(M)$, the most intrinsic of which is actually
$\chi(M)-(\chi(\Wbar)-\chi(C))=0$ (see below for the reason).

Now, given $\xi,\xi'\in\euls(M,\pp)$ we can choose generic representatives $v,v'$, so
that the set of points of $M$ where $v'=-v$ is a union of loops contained in the
interior of $M$. A standard procedure allows to give these loops a canonical orientation,
thus getting an element $\alphas(\xi,\xi')\in H_1(M;\mz)$. The following result is easily
obtained along the lines of the well-known analogue for closed manifolds.

\begin{lem}
$\alphas$ is well-defined and turns $\euls(M,\pp)$ into an affine space over $H_1(M;\mz)$.
\end{lem}

A (finite) cellularization $\cc$ of $M$ 
is called {\em suited} to $\pp$ if $V\cup C$ is
a subcomplex, so $W$ and $B$ are unions of cells. 
Here and in the sequel by ``cell'' we will always mean an {\em open} one.
Let such a $\cc$ be given.
For $\sigma\in\cc$ define $\index(\sigma)=(-1)^{{\rm dim}(\sigma)}$. We define
$\eulc(M,\pp)_\cc$ as the set of equivalence classes of integer singular 1-chains $z$ in $M$
such that 
$$\partial z=\sum_{\sigma\subset M\setminus (W\cup V)}\index(\sigma)\cdot p_\sigma$$
where $p_\sigma\in\sigma$ for all $\sigma$. Two chains $z$ and $z'$ 
with $\partial z=\sum\index(\sigma)\cdot p_\sigma$ and 
$\partial z'=\sum\index(\sigma)\cdot p'_\sigma$
are defined to be equivalent if there exist $\delta_\sigma:([0,1],0,1)\to
(\sigma,p_\sigma,p'_\sigma)$ such that
$$z-z'+\sum_{\sigma\subset M\setminus (W\cup V)}\index(\sigma)\cdot \delta_\sigma$$
represents $0$ in $H_1(M;\mz)$. Elements of
$\eulc(M,\pp)_\cc$ are called {\em combinatorial Euler structures} relative to $\pp$ and $\cc$,
and their representatives are called {\em Euler chains}.
The definition implies that, for $\xi,\xi'\in\eulc(M,\pp)_\cc$, their difference $\xi-\xi'$
can be defined as an element $\alphac(\xi,\xi')$ of $H_1(M;\mz)$. The following is easy:

\begin{lem}
$\eulc(M,\pp)_\cc$ is non-empty if and only if $\chi(\Wbar)=\chi(M)$, and in this case
$\alphac$ turns it into an affine space over $H_1(M;\mz)$.
\end{lem}

\noindent Since $\Wbar=W\cup V\cup C$, the alternating sum of dimensions of cells in $W\cup V$ is 
intrinsically interpreted as $\chi(\Wbar)-\chi(C)$, which explains why the most meaningful
way to write the relation $\chi(\Wbar)=\chi(M)$ is $\chi(M)-(\chi(\Wbar)-\chi(C))=0$.
From now on we will always assume that this relation holds.
Turaev~\cite{turaev:Euler} only considers the case where $V=C=\emptyset$, so 
$W=\overline{W}$ and $B=\overline{B}$, and our relation takes the usual form
$\chi(M,W)=0$.
The following result was established by Turaev in~\cite{turaev:Euler}
in his setting, but the proof
extends {\em verbatim} to our context, so we omit it.
Only the first assertion is hard.
We state the other two because we will use them.

\begin{prop}\label{combin:prop:statement}
\begin{enumerate}
\item If $\cc'$ is a subdivision of $\cc$ then there exists a canonical
$H_1(M;\mz)$-isomorphism $\eulc(M,\pp)_\cc\to\eulc(M,\pp)_{\cc'}$. In particular
$\eulc(M;\mz)$ is canonically defined up to $H_1(M;\mz)$-isomorphism independently 
of the cellularization.
\item\label{connected:spider:point} If 
$\cc$ is a cellularization of $M$ suited to $\pp$
and $x_0\in M$ is an assigned point, any element of $\eulc(M,\pp)$ can be represented,
with respect to $\cc$, as a sum $\sum_{\sigma\subset M\setminus(W\cup V)}
\index(\sigma)\cdot \beta_\sigma$
with $\beta_\sigma:([0,1],0,1)\to(M,x_0,\sigma)$.
\item\label{bary:sub:point} If $\calt$ is a triangulation of $M$ suited to $\pp$,
any element of $\eulc(M,\pp)$ can be represented, with respect to $\calt$, as 
a simplicial $1$-chain in the first barycentric subdivision of $\calt$.
\end{enumerate}
\end{prop}

\noindent Our first main result, proved in Section~\ref{proofs},
is the extension to the case under consideration of
Turaev's correspondence between $\eulc$ and $\euls$.

\begin{teo}\label{reconstruction:statement}
There exists a canonical $H_1(M;\mz)$-equivariant isomorphism 
$$\Psi:\eulc(M,\pp)\to\euls(M,\pp).$$
\end{teo}

\noindent The definition of $\Psi$ is based on an explicit geometric construction,
but its bijectivity is only established through $H_1(M;\mz)$-equivariance.
As already mentioned in the introduction, this makes in general a very difficult task to
determine the inverse of $\Psi$. One of the features of this paper is the description
of $\Psi^{-1}$ in terms of the combinatorial encoding of fields by 
means of branched spines: Theorem~\ref{spider:structure:teo} 
describes $\Psi^{-1}$ when $\pp$ is concave,
and Theorem~\ref{diagram:commutes} shows that from a general $\pp$ we can
effectively pass to a unique convex $\pp$, and hence to a unique concave $\pp$,
and conversely.

In view of Theorem~\ref{reconstruction:statement}, when no confusion risks to arise,
we shortly write $\eul(M,\pp)$ for either $\euls(M,\pp)$ or $\eulc(M,\pp)$,
and $\alpha$ for the map giving the affine $H_1(M;\mz)$-structure on this space.
Before turning to torsions, as announced in the introduction, we show
that (pseudo-)Legendrian links naturally define Euler structures of the type we are considering.

\begin{rem}\label{links:emerge}
{\em Assume $M$ is closed, let $v$ be a field on $M$ and let $L$ be a link in $M$
transverse to $v$. If we take a small enough regular neighbourhood $U(L)$ of $L$, 
the field $v$ will be tangent to $\partial U(L)$ only along two lines on each component,
and the tangency, viewed from the exterior $E(L)=M\setminus U(L)$, has concave type
(see the cross-section in Fig.~\ref{link:ext}).
\begin{figure}
\centerline{\psfig{file=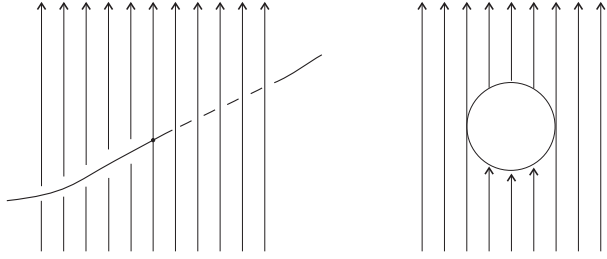,width=8cm}}
\caption{\label{link:ext} Concave tangency on a link exterior.}
\end{figure}
This shows that the triple $(M,v,L)$ defines an element of $\xi(M,v,L)$ of $\euls(E(L),\pp)$, 
where $\pp=(W,B,\emptyset,C)$
depends on the framing induced by $v$ on $L$. Note that if $\eta$ 
is a cooriented plane distribution and $L$ is tangent to $\eta$ then the definition
$\xi(M,\eta^\perp,L)$ applies. When $\eta$ is a contact structure $L$ is called a Legendrian link, so we will call it {\em pseudo-Legendrian} in general. In Section~\ref{link:section}
we shall see that $\xi(M,\eta^\perp,L)$ can be used to construct non-trivial invariants for
pseudo-Legendrian isotopy classes of knots.}
\end{rem}

\begin{figure}
\centerline{\psfig{file=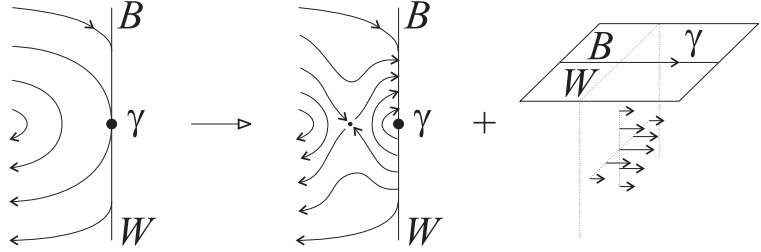,width=10cm}}
\caption{\label{conc:to:conv} Turning a concave tangency circle $\gamma$ into a convex one: the 
apparent singularity in the cross-section is removed by adding a small bell-shaped field 
directed parallel to $\gamma$, {\em i.e.}~orthogonal to the cross-section.}
\end{figure}

\paragraph{Convex Euler structure associated to an arbitrary one.}
Let $M$ and $\pp=(W,B,V,C)$ be as in the definition of $\eul(M,\pp)$.
The pattern $\theta(\pp)=(W,B,V\cup C,\emptyset)$ is a convex one canonically
associated to $\pp$. We define a map 
$$\Thetas:\euls(M,\pp)\to\euls(M,\theta(\pp))$$
as geometrically described in Fig.~\ref{conc:to:conv}. Concerning
this figure, note that the loops in $C$ can be oriented as 
components of the boundary of $B$, which is oriented as a subset of the boundary of $M$.

\begin{lem}\label{concave:convex:smooth}
$\Thetas$ is a well-defined $H_1(M;\mz)$-equivariant bijection.
\end{lem}

\dim{concave:convex:smooth}
The first two properties are easy and imply the third property. The inverse of $\Thetas$
may actually be described geometrically by a figure similar to Fig.~\ref{conc:to:conv}, but we leave this to the reader.
\finedim{concave:convex:smooth}
 
We define now a combinatorial version of $\Thetas$. Consider a cellularization $\cc$ suited to 
$\pp$, and denote by $\gamma_1,\dots,\gamma_n$ the 1-cells contained in $C$. We choose the 
parameterizations $\gamma_j:(0,1)\to C$ so that they respect the natural orientation of $C$ 
already discussed above, and we extend the $\gamma_j$ to $[0,1]$, without changing notation.
Now let $z$ be an Euler chain relative to $\pp$. It 
easily seen that $z-\sum_{j=1}^n\gamma_j\ristr{[1/2,1]}$ is an Euler chain
relative to $\theta(\pp)$. Setting 
$$\Thetac([z])=\left[z-\sum_{j=1}^n\gamma_j\ristr{[1/2,1]}\right]$$ 
we get a map $\Thetac:\eulc(M,\pp)\to\eulc(M,\theta(\pp))$.

\begin{lem}\label{concave:convex:comb}
$\Thetac$ is a well-defined $H_1(M;\mz)$-equivariant bijection.
\end{lem}

\dim{concave:convex:comb}
Again, the first two properties are easy and imply the third one. 
\finedim{concave:convex:comb}

In Section~\ref{proofs} we will see the following:

\begin{teo}\label{diagram:commutes}
If $\Psi$ is the reconstruction map of Theorem~\ref{reconstruction:statement} then the
following diagram is commutative:
$$\matrix{
\eulc(M,\pp) & \stackrel{\Thetac}{\longrightarrow} & \eulc(M,\theta(\pp))\phantom{.} \cr
\Psi\downarrow\phantom{\Psi} & & \phantom{\Psi}\downarrow\Psi \cr
\euls(M,\pp) & \stackrel{\Thetas}{\longrightarrow} & \euls(M,\theta(\pp)). \cr}$$
\end{teo}

Using this result we will sometimes just write 
$\Theta:\eul(M,\pp)\to\eul(M,\theta(\pp))$.

\paragraph{Torsion of a (convex) Euler structure.}
Let us first briefly review the algebraic setting~\cite{milnor} in which 
torsions can be defined. 
We consider a ring $\Lambda$ with unit, with the property
that if $n$ and $m$ are distinct positive integers then $\Lambda^n$ and 
$\Lambda^m$ are not isomorphic as $\Lambda$-modules. The
Whitehead group $K_1(\Lambda)$ is defined as the
Abelianization of $\GL_\infty(\Lambda)$, and $\Kbar_1(\Lambda)$ is the
quotient of $K_1(\Lambda)$ under the action of $-1\in\GL_1(\Lambda)=
\Lambda_*$. Now consider a convex boundary pattern
$\pp=(W,B,V,\emptyset)$ on a manifold $M$, take
a representation $\varphi:\pi_1(M)\to\Lambda_*$, and consider the $\Lambda$-modules $H^\varphi_i(M,W\cup V)$ of relative twisted homology 
(see Section~\ref{torsion:def:section} for a reminder on the definition). 
Notice that $W\cup V=\Wbar$, so if we have a cellularization of $M$ suited to $\pp$
then $W\cup V$ is a (closed) subcomplex, and we can use
the cellular theory to compute $H^\varphi_i(M,W\cup V)$. This is the reason for
considering {\em convex} boundary patterns.

Assume now that 
$H^\varphi_i(M,W\cup V)$ is free, and choose a $\Lambda$-basis $\hbasis_i$. Then a torsion
$\tau^\varphi(M,W\cup V,\hbasis)\in\Kbar_1(\Lambda)/\varphi(\pi_1(M))$ can be defined as in~\cite{milnor}. Here the action of $\varphi(\pi_1(M))$ has to be taken into account
because of the ambiguity of the choice of liftings to the universal cover
of the cells of $M\setminus(W\cup V)$. It was an intuition of Turaev~\cite{turaev:Euler}, which
we extend in this paper to include the case of simple tangency, 
that an Euler structure $\xi\in\eul(M,\pp)$	
can be used to get rid of the action of $\varphi(\pi_1(M))$. More precisely we will
show below the following:

\begin{teo}\label{torsion:existence:statement}
In the above situation a torsion $\tau^\varphi(M,\pp,\xi,\hbasis)$ can be defined. 
The reduction modulo $\varphi(\pi_1(M))$ of $\tau^\varphi(M,\pp,\xi,\hbasis)$
gives $\tau^\varphi(M,W\cup V,\hbasis)$. Moreover if $\xi,\xi'\in\eul(M,\pp)$ then
\begin{equation}
\tau^\varphi(M,\pp,\xi',\hbasis)=\tau^\varphi(M,\pp,\xi,\hbasis)\cdot 
\overline{\varphi}(\alpha(\xi',\xi)).\label{equivariance:formula}
\end{equation}
\end{teo}

For a formal definition of $\overline{\varphi}:H_1(M;\mz)\to\Kbar_1(\Lambda)$ see
Section~\ref{torsion:def:section}. A self-contained definition of 
$\tau^\varphi(M,\pp,\xi,\hbasis)$ will also be given in Section~\ref{torsion:def:section}.
For the readers already acquainted with~\cite{milnor}, we mention the key point:
given a cellularization of $M$ suited to $\pp$,
a preferred family of liftings for the cells in $M\setminus(W\cup V)$ is found
by representing $\xi$ by a ``connected spider'' as in
point~\ref{connected:spider:point} of Proposition~\ref{combin:prop:statement},
lifting the spider starting from an arbitrary
lifting of its head, and defining the preferred cell-liftings 
as those containing the ends of the legs of the lifted spider.

Theorem~\ref{torsion:existence:statement} only applies to convex patterns,
but if $\pp$ is not convex we can use the canonical map 
$\Theta:\eul(M,\pp)\to\eul(M,\theta(\pp)$ and define  
$$\tau^\varphi(M,\pp,\xi,\hbasis)=\tau^\varphi(M,\theta(\pp),\Theta(\xi),\hbasis).$$
Of course the equivariance formula~(\ref{equivariance:formula}) still holds.

One of the important features of
Theorem~\ref{torsion:existence:statement} is that if we start from a {\em combinatorial}
representative of $\xi$ then the computation of $\tau^\varphi(M,\pp,\xi,\hbasis)$
is (in principle) algorithmic, provided we start from an explicit description
of the universal cover of $M$ (or the maximal Abelian cover, which is often easier,
when $\Lambda$ is commutative). 

The next result follows directly from the definition but is nonetheless worth stating, because it shows how torsions may be used to distinguish triples $(M,\pp,\xi)$ from each other (see
Section~\ref{link:section} for a relevant consequence).

\begin{prop}\label{homeo:action}
Let $f: M\to M'$ be a homeomorphism, consider $\xi\in\eul(M,\pp)$, 
$\varphi:\pi_1(M)\to\Lambda_*$ and a $\Lambda$-basis $\hbasis$ of 
$H^\varphi_*(M,\overline{W})$. Then
$$\tau^{\varphi\circ f_*^{-1}}(M',f_*(\pp),f_*(\xi),f_*(\hbasis))=
\tau^\varphi(M,\pp,\xi,\hbasis).$$
\end{prop} 

\section{Torsion of pseudo-Legendrian knots}\label{link:section}

We fix in this section a compact oriented manifold $M$ and a boundary pattern $\pp$ on $M$.
The boundary of $M$ may be empty or not. Recall that if $v$ is a vector field on $M$ and 
$K$ is a knot in $\interior(M)$, we have defined $K$ to be pseudo-Legendrian in $(M,v)$ 
if $v$ is transversal to $K$. We will also call $(v,K)$ a pseudo-Legendrian pair.
Having fixed $\pp$, we will only consider fields $v$ 
compatible with $\pp$. The aim of this section is to show how torsions can be applied to
distinguish pseudo-Legendrian knots. Some of the results we will establish hold
also for links, but we will stick to knots for the sake of 
simplicity. First, we need to spell out the equivalence relation
which we consider.

Let $v_0,v_1$ be compatible with $\pp$ and let $K_0,K_1$ be pseudo-Legendrian in
$(M,v_0)$ and $(M,v_1)$ respectively. We define $(v_0,K_0)$ to be {\em weakly equivalent}
to $(v_1,K_1)$ if there exist a homotopy $(v_t)_{t\in[0,1]}$ through fields compatible
with $\pp$ and an isotopy $(K_t)_{t\in[0,1]}$ such that $K_t$ is transversal
to $v_t$ for all $t$. If $v_0=v_1$ then $K_0$ and $K_1$ are
called {\em strongly equivalent} if the homotopy $(v_t)$ can be chosen to be constant.

\begin{rem}{\em
Of course strong equivalence implies weak equivalence. Weak equivalence is the natural
relation to consider on pseudo-Legendrian pairs $(v,K)$, while strong equivalence
is natural for pseudo-Legendrian knots in a fixed $(M,v)$. We will see
that torsion provides obstructions to weak (and hence to strong) equivalence.}
\end{rem}

Now let $K$ be pseudo-Legendrian in $(M,v)$ and note that $v$ turns 
$K$ into a framed knot, which we will denote by $K^{(v)}$.
The framed-isotopy class of $K^{(v)}$ is of course invariant under weak
equivalence, so we will only try to distinguish knots which are framed-isotopic 
to each other. As already mentioned, the idea is to restrict $v$ to the exterior of
$K$ and consider the induced Euler structure. A technical
subtlety arises here, because the comparison class of two such Euler structures
coming from framed-isotopic knots cannot be computed directly. It will
turn out that the action of a group must be taken into account.
However, we will see that for important classes
of knots this action can actually be neglected.

\paragraph{Euler structures on knot exteriors.}
For a knot $K$ in $M$ we consider a (closed) tubular neighbourhood $U(K)$ of $K$
in $M$ and we define $E(K)$ as the closure of the complement of $U(K)$.
If $F$ is a framing on $K$ we extend the boundary pattern $\pp$ previously fixed
on $M$ to a boundary pattern $\pp(K^F)$ on $E(K)$, by splitting $\partial U(K)$
into a white and a black longitudinal annuli, the longitude being the one defined by
the framing $F$. As a direct application of Proposition~\ref{p:h:formula} one sees
that $\eul(E(K),\pp(K^F))$ is non-empty (assuming $\eul(M,\pp)$ to be non-empty).

A convenient way to think of $\pp(K^F)$ is as follows. The framing $F$ 
determines a transversal vector field along $K$. If we extend this field near
$K$ and choose  $U(K)$ small enough then the pattern we see on $\partial U(K)$
is exactly as required. With this picture in mind, it is clear that if $K$ is
pseudo-Legendrian in $(M,v)$, where $v$ is compatible with $\pp$, then the
restriction of $v$ to $E(K)$ defines an element
$$\xi(v,K)\in\eul(E(K),\pp(K^{(v)}).$$ 
(This notation is consistent with that
previously used, because in this section we are considering $M$ to be fixed.)

\paragraph{Group action on Euler structures.}
Consider a knot $K$ and a
self-diffeomorphism $f$ of $E(K)$ which is the identity near $\partial E(K)$.
Then $f$ extends to a self-diffeomorphism $\widehat f$ of $M$, where 
$\widehat f\ristr{U(K)}={\rm id}_{U(K)}$. We define $G(K)$ as the group of all such $f$'s
with the property that $\widehat f$ is isotopic to the identity on $M$.
Elements of $G(K)$ are regarded up to isotopy relative to $\partial E(K)$.
If $F$ is a framing on $K$ then the pull-forward of vector fields induces
an action of $G(K)$ on $\eul(E(K),\pp(K^{(v)})$. We will now see 
that an obstruction to weak equivalence can be expressed in terms
this group action.

Let $(v_0,K_0)$ and $(v_1,K_1)$ be pseudo-Legendrian pairs in $M$, and assume
that  $K_0^{(v_0)}$ is framed-isotopic to $K_1^{(v_1)}$ under a diffeomorphism
$f$ relative to $\partial M$. Using the restriction of $f$ and the pull-back of
vector fields we get a bijection
$$f^*:\eul(E(K_1),\pp(K_1^{(v_1)}))\to\eul(E(K_0),\pp(K_0^{(v_0)})).$$

\begin{prop}\label{group:obstr}
Under the current assumptions, if $(v_0,K_0)$ and $(v_1,K_1)$ are weakly equivalent to each
other then $f^*(\xi(v_1,K_1))$ belongs to the $G(K_0)$-orbit of $\xi(v_0,K_0)$ in
$\eul(E(K_0),\pp(K_0^{(v_0)})$.
\end{prop}

\dim{group:obstr}
By assumption $K_0,K_1$ and $v_0,v_1$ embed in continuous families $(K_t)_{t\in[0,1]}$
and $(v_t)_{t\in[0,1]}$, where $v_t$ is transversal to $K_t$ for all $t$. Now
$(K_t^{(v_t)})_{t\in[0,1]}$ is a framed-isotopy, so 
there exists a continuous family $(g_t)_{t\in[0,1]}$ of diffeomorphisms of $M$
fixed on $\partial M$ and such that $g_0={\rm id}_M$ and $g_t(K_0^{(v_0)})=K_t^{(v_t)}$. 
So we get a map
$$[0,1]\ni t\mapsto
\alpha(\xi(v_0,K_0),g_t^*(\xi(v_t,K_t)))\in H_1(E(K_0);\mz).$$
Since $H_1(E(K_0);\mz)$ is discrete and the map is continuous,
we deduce that the map is identically 0. So $g_1^*(\xi(v_1,K_1))=\xi(v_0,K_0)$.
Now
$$f^*(\xi(v_1,K_1))=(f^*\compo (g_1)_*\compo g_1^*)(\xi(v_1,K_1))=
(f^{-1}\compo g_1)_*(\xi(v_0,K_0))$$
and the conclusion follows because $f^{-1}\compo g_1$ defines an element of $G(K_0)$.
\finedim{group:obstr}

The group $G(K)$ is in general rather difficult to understand (see~\cite{hatcher}), 
so we introduce a special terminology for the case where its action
can be neglected. We will say that a framed knot $K^F$ is {\em good} if
$G(K)$ acts trivially on $\eul(E(K),\pp(K^F))$. If $K^F$ is good for all framings $F$,
we will say that $K$ itself is good.
The following are easy examples of good knots:
\begin{itemize} 
\item $M$ is $S^3$ and $K$ is the trivial knot;
\item $M$ is a lens space $L(p,q)$ and $K$ is the core of one of the handlebodies
of a genus-one Heegaard splitting of $M$.
\end{itemize} 
\noindent The reason is that in both cases $E(K)$ is a solid torus, and we know that an automorphism of
the solid torus which is the identity on the boundary is isotopic 
to the identity relatively to the boundary, so $G(K)$ is trivial. The next three
results show that on one hand $G(K)$ is very seldom trivial, but on the other hand 
many knots are good. We will give proofs in the sequel, after introducing
some extra notation. In the statements, by `$E(K)$ is hyperbolic' 
we mean `$\interior(E(K))$ is complete, finite-volume hyperbolic.'

\begin{prop}\label{seldom:triv}
If $M$ is closed and $E(K)$ is hyperbolic then $G(K)$ is non-trivial.
\end{prop}

\begin{teo}\label{many:good}
If $M$ is closed, $E(K)$ is hyperbolic and 
either ${\rm Out}(\pi_1(E(K)))$ is trivial or $H_1(E(K);\mz)$ is torsion-free
then $K$ is good.
\end{teo}

\begin{teo}\label{homology:good}
If $M$ is a homology sphere then every knot in $M$ is good.
\end{teo}

The next result, which follows directly from Proposition~\ref{group:obstr},
the definition of goodness, and Proposition~\ref{homeo:action},
shows that for good knots torsion can be used as an
obstruction to weak (and hence strong) equivalence.

\begin{prop}\label{torsion:obstr}
Let $(v_0,K_0)$ and $(v_1,K_1)$ be pseudo-Legendrian pairs in $M$, and assume that 
$K_0^{(v_0)}$ is framed-isotopic to $K_1^{(v_1)}$ under a diffeomorphism $f$ relative
to $\partial M$. Suppose that $K_0^{(v_0)}$ is good, and that for some representation
$\varphi:\pi_1(E(K_0))\to\Lambda$ and some $\Lambda$-basis $\hbasis$ of
$H^\varphi_*(E(K_0),\overline{W(\pp(K_0^{(v_0)}))})$ we have
\begin{equation}\label{tors:obs:eq}
\tau^\varphi(E(K_0),\pp(K_0^{(v_0)}),\xi(v_0,K_0),\hbasis)\neq
\tau^{\varphi\compo f_*^{-1}}(E(K_1),\pp(K_1^{(v_1)}),\xi(v_1,K_1),f_*(\hbasis)).
\end{equation}
Then $(v_0,K_0)$ and $(v_1,K_1)$ are not weakly equivalent.
\end{prop}

\begin{rem}{\em
\begin{enumerate}
\item The right-hand side of equation~(\ref{tors:obs:eq}) actually equals
$$\tau^\varphi(E(K_0),\pp(K_0^{(v_0)}),f^*(\xi(v_1,K_1)),\hbasis),$$
but we have written it as it stands in order to use only the action of $f$
on the fundamental group and on the twisted homology, not on Euler chains.
Using the technology described in Section~\ref{spines:section}, both sides
of the equation can be computed in practice, at least when $\Lambda$ is commutative.
\item An obstruction in terms of torsion may be given also for non-good knots,
but the statement would become awkward and nearly impossible to apply, so we have 
refrained from giving it.
\item If equation~(\ref{tors:obs:eq}) holds for some basis $\hbasis$ then it holds
for any basis.
\end{enumerate}}
\end{rem}

\noindent To conclude this paragraph we note that using the technology of
Turaev~\cite{turaev:Euler},  one can actually see that the action on Euler
structures of an automorphism  is invariant under {\em homotopy} (not only
isotopy) relative to the boundary. We will not use this fact.

\paragraph{Good knots.} We introduce now some notation needed for the
proofs of Proposition~\ref{seldom:triv} and Theorem~\ref{many:good} (for 
Theorem~\ref{homology:good} we will use a different approach, see below).
Recall that $(M,\pp)$ is fixed for the whole section. We temporarily fix also
a framed knot $K^F$ in $M$, a regular neighbourhood $U$ of $K$, and we denote by $T$ the
boundary torus of $U$. On $T$ we consider 1-periodic coordinates $(x,y)$
such that $x\mapsto (x,0)$ is a meridian of $U$ and $y\mapsto (0,y)$ is 
a longitude compatible with $F$. We denote a collar of $T$ in $E(K)$ by 
$V$ and parameterize $V$ as $T\times[0,1]$, where $T=T\times\{0\}$. We consider on
$[0,1]$ a coordinate $s$. For $p,q\in\mz$ we define automorphisms
$\dd_{(p,q)}$ of $E(K)$ as follows. 
Each $\dd_{(p,q)}$ is supported in $V$, and on $V$, using the
coordinates just described, it is given by 
$$\dd_{(p,q)}(x,y,s)=(x+p\cdot s, y +q\cdot s, s).$$
We will call such a map a {\em Dehn twist}. It is easy to verify that
the extension of $\dd_{(p,q)}$ to $M$ is isotopic to the identity of $M$.
Note that $\dd_{(p,q)}$ is actually not smooth on $T\times\{1\}$, but 
we can consider some smoothing and identify $\dd_{(p,q)}$
to an element of $G(K)$, because the equivalence class is
independent of the smoothing.

\dim{seldom:triv}
We show that $\dd_{(p,q)}$ is non-trivial in $G(K)$ for all $(p,q)\neq(0,0)$. 
Fix the basepoint $a_0=(0,0)\in T$ for the fundamental groups of $T$ and $E(K)$. Then 
$\dd_{(p,q)}$ acts on $\pi_1(E(K),a_0)$ as the conjugation by $i_*(p,q)$, where 
$i:T\to E(K)$ is the inclusion and $(p,q)\in\mz\times\mz=\pi_1(T,a_0)$.
If $\dd_{(p,q)}$ is trivial in $G(K)$, {\em i.e.} it is isotopic to the identity
relatively to $\partial E(K)$, in particular it acts trivially on 
$\pi_1(E(K),a_0)$. This implies that $i_*(p,q)$ is in the centre of $\pi_1(E(K),a_0)$.
Now it follows from hyperbolicity that this centre is trivial and $i_*$ is injective,
whence the conclusion.
\finedim{seldom:triv}

\noindent The proof of Theorem~\ref{many:good} will rely on properties of hyperbolic
manifolds and on the following fact, which we consider to be quite remarkable
(note that the 2-dimensional analogue,
which may be stated quite easily, is false).
Remark that the result applies in particular to Dehn twists.

\begin{prop}\label{collar:trivial}
If $[f]\in G(K)$ and $f$ is supported in the collar $V$ of $\partial U$ then $[f]$
acts trivially on $\eul(E(K),\pp(K^F))$.
\end{prop}

\dim{collar:trivial}
Consider a vector field $v$ on $E(K)$ compatible with $\pp(K^F)$. Since $v$ and $f_*(v)$
differ only on $V$, their difference belongs to the image of $H_1(V;\mz)$ in
$H_1(E(K);\mz)$. So we may as well assume that $E(K)=V$, {\em i.e.} $M$ is
the solid torus $U\cup V$.

By contradiction, let $\xi\in\eul(V,\pp(K^F))$ be such that 
$\alpha(\xi,(\dd_{(p,q)})_*(\xi))$ is non-zero in $H_1(V;\mz)$, so it is given
by $k\cdot[\gamma]$ for some $k\in\mz\setminus\{0\}$ and some simple closed
curve $\gamma$ on $T\times\{1\}\subset\partial V$. Let us now take another
simple closed curve $\delta$ on $T\times\{1\}$ which intersects $\gamma$
transversely at one point. Let us define $N$ as the manifold obtained by
attaching the solid torus to $V$ along $T\times\{1\}$, in such a way that the
meridian of the solid torus gets identified  with $\delta$. Note that $N$ is
again a solid torus and that the homology class of $\gamma$ in
$H_1(N;\mz)\cong\mz$ is a generator. Now we can apply
Proposition~\ref{p:h:formula} to extend $\xi$ to an Euler structure $\xi_N$ on
$N$. Moreover we can extend $f$ to an automorphism $g$ of $N$ which is the
identity on $\partial N=T\times\{0\}$. Now by construction
$\alpha(\xi_N,g_*(\xi_N))$ equals $k\cdot[\gamma]$ in $H_1(N;\mz)\cong\mz$, so
it is non-zero. But $g$ is isotopic to the identity of $N$ relatively to the
boundary of $N$, so we have a contradiction. \finedim{collar:trivial}

For the proof of Theorem~\ref{many:good} we will also need the
following easy fact. 

\begin{lem}\label{power:action}
Let $f$ be an automorphism of $M$ relative to $\partial M$, and consider the
induced automorphisms of $H_1(M;\mz)$ and $\eul(M,\pp)$, both denoted by $f_*$. Then:
$$\alpha(f_*(\xi_0),f_*(\xi_1))=f_*(\alpha(\xi_0,\xi_1)),\qquad
\forall\xi_0,\xi_1\in\eul(M,\pp).$$
\end{lem}

\dim{power:action}
Take representatives of $\xi_0$ and $\xi_1$ such that $\alpha(\xi_0,\xi_1)$ can be
viewed as the anti-parallelism locus. The formula is then obvious.
\finedim{power:action}

\dim{many:good} Consider $[f]\in G(K)$. It follows from the work of Johansson
(see~\cite{hatcher}) that, under the assumption that $E(K)$ is hyperbolic, the
group generated by Dehn twists has finite index in the mapping class group of
$E(K)$ relative to the boundary. More precisely, the quotient group can be
identified to a subgroup of ${\rm Out}(\pi_1(E(K))$, which is finite as a
consequence of Mostow's rigidity. If ${\rm Out}(\pi_1(E(K))$ is trivial then
$[f]$ is equivalent to a Dehn twist, so $f$ acts trivially on
$\eul(E(K),\pp(K^F))$ by Proposition~\ref{collar:trivial}. 

We are left to deal with the case where $H_1(E(K);\mz)$ is torsion-free. By
Johansson's result, there exists an integer $n$ such that $f^n$ acts trivially
on $\eul(E(K),\pp(K^F))$. Consider now $\xi\in\eul(E(K),\pp(K^F))$, and set
$\alpha=\alpha(\xi,f_*(\xi))$.  We must show that  $\alpha=0$. We denote by
$\widehat\alpha$ the image of $\alpha$ in $H_1(M;\mz)$, and by  $\widehat f$
the extension of $f$ to $M$. Since $\widehat f$ is isotopic to the identity,
we  have $\widehat f_*(\widehat\alpha)=\widehat\alpha$. If we take an oriented
1-manifold $a$ representing $\alpha$ and disjoint from $\partial U(K)$, this
means that there exists an  oriented surface $\Sigma$ in $M$ such that
$\partial\Sigma=a\cup(-f(a))$. Up to isotopy we can  assume that $\Sigma$
intersects $\partial U(K)$ transversely in a union of circles. This shows that
$f_*(\alpha)=\alpha+k\cdot\mu$, where $\mu$ is the meridian of $T$. Note that
$f_*(\mu)=\mu$, so for all integers $m$ we have $f_*^m(\alpha)=\alpha+m\cdot
k\cdot\mu$. Now, using Lemma~\ref{power:action}, we have:
	\begin{eqnarray*}
0 & = & \alpha(\xi,f_*^n(\xi)=\sum_{m=0}^{n-1}\alpha(f_*^m(\xi),f_*^{m+1}(\xi))\\
& = & \sum_{m=0}^{n-1} f_*^m(\alpha(\xi,f_*(\xi)))=\sum_{m=0}^{n-1}f_*^m(\alpha)=
\sum_{m=0}^{n-1}(\alpha+m\cdot k\cdot\mu)\\
& = & n\cdot \alpha +{n(n-1)\over 2}\cdot k\cdot \mu.
\end{eqnarray*}
This shows that $2\cdot \alpha+(n-1)\cdot k\cdot \mu$ is a torsion element of $H_1(E(K);\mz)$,
so it is null by assumption. So $(1-n)\cdot k\cdot \mu=2\cdot \alpha$. If we apply $f_*$ to
both sides of this equality we get
$(1-n)\cdot k\cdot f_*(\mu)=2\cdot f_*(\alpha)$. Using the equality again and the relations
$f_*(\mu)=\mu$ and $f_*(\alpha)=\alpha+k\cdot\mu$ we get
$$(1-n)\cdot k\cdot \mu=2\cdot\alpha+2\cdot k\cdot\mu=(1-n)\cdot k\cdot
\mu+2\cdot k\cdot\mu.$$
Therefore $k\cdot\mu$ is a torsion element, and hence null. But
$2\cdot\alpha=(1-n)\cdot k\cdot\mu$, so also $\alpha$ is null.\finedim{many:good}

\paragraph{Torsion and rotation number, and more good knots.} 
We will show in this section that for a contact 
structure in a homology sphere the rotation number of a Legendrian knot can be expressed 
in terms of Euler structures on the complement. This will imply that torsion 
essentially contains the rotation number, and it will allow us to show
that in a homology sphere all knots are good (Theorem~\ref{homology:good}).

To begin, we note that the definition of the rotation number, classically defined in the
contact case, actually extends to the situation we are considering. Since we will need
this definition, we recall it. Let $M$ be a homology sphere, let
$v$ be a field on $M$ and let $K$ be an oriented pseudo-Legendrian knot 
in $(M,v)$. Take a plane field $\eta$ transversal to $v$ and tangent to $K$, 
and a Seifert surface
$S$ for $K$. Up to isotopy of $S$ we can assume that $\eta$ is tangent to
$S$ only at isolated points. Then ${\rm rot}_v(K)$ is the sum of a contribution
for each of these tangency points $p$. Define ${\rm o}(p)$ to be $+1$ if
$\eta_p=T_pS$ and $-1$ if $\eta_p=-T_pS$. If $p\in\partial S=K$ then $p$ contributes
just with ${\rm o}(p)$. If $p\in\interior(S)$ we can consider near $p$ a section of 
$\eta\cap TS$ which vanishes at $p$ only, and denote by ${\rm i}(p)$ its index. Then
$p$ contributes to ${\rm rot}_v(K)$ with ${\rm o}(p)\cdot{\rm i}(p)$.

It is quite easy to see that the resulting number is indeed independent from $\eta$
and $S$. Moreover ${\rm rot}_v(K)$ is unchanged under homotopies of $v$ relative to $K$,
and local modifications away from $K$, so we can actually define
${\rm rot}_\xi(K)$ where $\xi=\xi(v,K)\in\eul(E(K),\pp(K^{(v)})$.

\begin{prop}\label{alpha:rot}
Let $M$ be a homology sphere, let $v$ be a field on $M$ 
and let $K_0$ and $K_1$ be oriented pseudo-Legendrian knots in $(M,v)$. 
Assume that there exists a framed-isotopy $f$ which maps $K_1^{(v)}$ to 
$K_0^{(v)}$. Identify $H_1(E(K_0);\mz)$ to $\mz$ using a meridian. Then:
$${\rm rot}_v(K_1)={\rm rot}_v(K_0)+2\alpha(f_*(\xi(v,K_1)),\xi(v,K_0)).$$
\end{prop}

\dim{alpha:rot}  Let $K:=K_0$, $v_0:=v$ and $v_1:=f_*(v)$. Note that $v_0$ and
$v_1$ coincide along $K$.  Of course ${\rm rot}_v(K_1)={\rm rot}_{v_1}(K)$. We
are left to show that 
$${\rm rot}_{v_1}(K)={\rm rot}_{v_0}(K)+2\alpha(\xi(v_1,K)),\xi(v_0,K)).$$ 
We can now homotope $v_0$ and $v_1$ away from $K$ until they differ only in
the  neighbourhood $W(L)$ of an oriented link $L$, and within this
neighbourhood they differ exactly by a ``Pontrjagin move'', as defined for
instance in~\cite{lnm}. Namely, $v_0$  runs parallel to $L$ in $W(L)$, while
$v_1$ runs opposite to $L$ on $L$ and has non-positive radial component on
$W(L)$ (see below for a picture). Note that $L$ represents
$\alpha(\xi(v_1,K)),\xi(v_0,K))$.

Let us choose now a Seifert surface $S$ for $K$ and a Riemannian metric on $M$,
and define $\eta_i=v_i^\perp$, for $i=0,1$. Since
$\eta_0\ristr{K}=\eta_1\ristr{K}$, the contributions along $K$ to ${\rm
rot}_{v_0}(K)$ and ${\rm rot}_{v_1}(K)$ are the same. Up to isotoping $S$ we
may assume that $L$ is transversal but never orthogonal to $S$. At the points
where $\eta_0$ is tangent to $S$ also $\eta_1$ is tangent to $S$, and the
contributions are the same. So ${\rm rot}_{v_1}(K)-{\rm rot}_{v_0}(K)$ is given
by the sum of the contributions of the tangency points of $\eta_1$ to $S$
within $W(L)$. We will show that each point of $L\cap S$ gives rise to exactly
two tangency points, which both contribute with $+1$ or $-1$ according to the
sign of the intersection of $L$ and $S$ at the point. This will show that ${\rm
rot}_{v_1}(K)-{\rm rot}_{v_0}(K)$ is twice the algebraic intersection of $L$
and $S$. This algebraic intersection is exactly the value of
$[L]=\alpha(\xi(v_1,K)),\xi(v_0,K))$ as a multiple of $[m]$, so the local
analysis at $L\cap S$ will imply the desired conclusion.

\begin{figure}
\centerline{\psfig{file=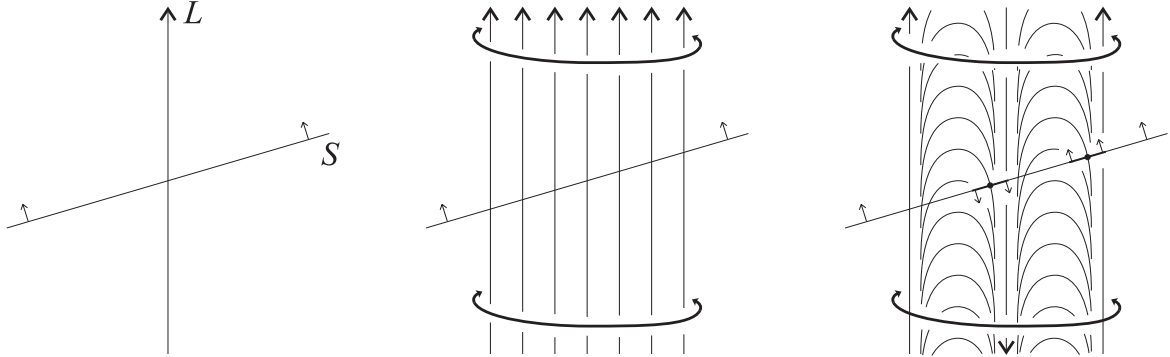,width=15.5cm}}
\caption{\label{rot:alpha:fig}Effect of the Pontrjagin move.}
\end{figure}

For the sake of simplicity we only examine a positive intersection point of $L$
and $S$. This is done in a cross-section in Fig.~\ref{rot:alpha:fig}, which
shows the local effect
of the move. The fields pictured both have a rotational symmetry, suggested in the
figure. The two tangency points which arise are a positive focus (on the right) and
a negative saddle (on the left), 
so the local contribution is indeed $+2$, and the proof is complete.
\finedim{alpha:rot}

\begin{rem}{\em The definition of rotation number and
Proposition~\ref{alpha:rot}
easily extend to the case of manifolds which are not homology spheres, 
by restricting to homologically trivial knots and choosing a relative homology class in the
complement.}
\end{rem}

\noindent We can now prove that in a homology sphere all knots are good.

\dim{homology:good}
Consider $[f]\in G(K)$, a framing $F$ on $K$ and $\xi\in\eul(E(K),\pp(K^F))$. We must
show that $f_*(\xi)=\xi$. Let $\xi=[v]$ and denote by $\widehat{v}$ the
obvious extension of $v$ to $M$. As above, let $\widehat{f}$ be the extension of $f$ to $M$.
During the proof of Proposition~\ref{alpha:rot} we have shown that
$${\rm rot}_{\widehat{f}_*(\widehat{v})}(K)-{\rm rot}_{\widehat{v}}(K)=
2\alpha(f_*(v),v).$$
But ${\rm rot}_{\widehat{f}_*(\widehat{v})}(K)$ is actually equal to 
${\rm rot}_{\widehat{v}}(K)$, because 
$\widehat{f}$ is the identity near $K$. Therefore
$f_*(v)$ and $v$ differ by a torsion element of
$H_1(E(K);\mz)\cong\mz$, so they are equal. By definition
$f_*(\xi)=[f_*(v)]$ and $\xi=[v]$, and the proof
is complete.
\finedim{homology:good}

Theorems~\ref{many:good} and~\ref{homology:good} provide a partial answer
to the problem of determining which knots are good. The general problem does not appear 
to be straight-forward, and we leave it for further investigation.
We will only show below an example of knot which is not good.

\paragraph{Knots distinguished by torsion.}
This paragraph is devoted to proving Theorem~\ref{informal:good:for:links}.
Within the proof we will need the following general
fact, which we state separately:

\begin{lem}\label{follow:homotopy}
Let $(v_t)_{t\in[0,1]}$ be a homotopy of non-singular vector fields on a 
$3$-manifold $M$, and let $K_0$ be a knot transversal to $v_0$. Then $K_0$ extends to
an isotopy $(K_t)_{t\in[0,1]}$ such that $K_t$ is transversal to $v_t$ for all $t$.
\end{lem}

\noindent This lemma can be established
using the classical methods of general position and obstruction theory,
and we leave it to the reader. We just mention that an easy alternative proof could also
be given in the framework of the theory of branched standard spines, using Theorem~\ref{closed:bijectivity:teo} and $\cont^1$ projections of knots (see~\cite{second:paper}). 

We state now  our main result, addressing the reader to~\cite{eliash} for the
definition of {\em overtwisted} contact structure. Before giving the proof
we discuss the consequences which are most relevant to us.

\begin{prop}\label{link:family:prop}
Let $(v,K)$ be a pseudo-Legendrian pair in $M$. For all
$$\gamma\in{\rm Ker}(i_*:H_1(E(K);\mz)\to H_1(M;\mz))$$
there exists a pseudo-Legendrian knot $K_\gamma$ in $(M,v)$ and an isotopy
$f:M\to M$ which maps $K_\gamma^{(v)}$ to $K^{(v)}$ such that 
$$\alpha(\xi(v,K),f_*(\xi(v,K_\gamma)))=\gamma.$$
Moreover, if $v$ is transversal to an assigned overtwisted contact structure $\eta$
and $K$ is Legendrian in $\eta$ then also $K_\gamma$ can be chosen to be
Legendrian in $\eta$.
\end{prop}

\begin{rem}\label{tors:really:OK}{\em  If $K^{(v)}$ is good and $\gamma\neq 0$,
the pairs $(v,K)$ and $(v,K_\gamma)$ are not weakly equivalent, and a torsion
tells them apart, in the sense that Proposition~\ref{torsion:obstr} applies.
(To see this, choose $\varphi:H_1(E(K);\mz)\to\Lambda_*$ with
$\varphi(\gamma)\neq1$, consider the induced representation of $\pi_1(E(K))$,
and apply formula~(\ref{equivariance:formula}) of
Theorem~\ref{torsion:existence:statement}.)} \end{rem}

\begin{rem}\label{tors:really:super} {\em When $M$ is a homology sphere, so
that $K$ is automatically good, the family of knots obtained from
Proposition~\ref{link:family:prop} is  parameterized by $\mz$, and we can choose
a representation $\varphi:\pi_1(E(K))\to\Lambda_*$  such that $\tau^\varphi$
takes a different value on each knot of the family. This shows in particular
that the knots are pairwise weakly inequivalent. In the contact case,  the
knots are pairwise framed-isotopic but not isotopic through Legendrian knots.
(For a proof, choose $\varphi$ such that $\varphi(1)$ has infinite order.)}
\end{rem}

\dim{link:family:prop}
We start by modifying
the field $v$ on $E(K)$ to a field $w$, without modification near $\partial E(K)$,
in such a way that $\alpha(v\ristr{E(K)},w)=\gamma\in H_1(E(K);\mz)$. 
This can be achieved by a ``Pontrjagin move'', as already used within 
the proof of~\ref{alpha:rot}. Let us spell out the steps to be followed:
\begin{enumerate}
\item Select an oriented link $L$ in the interior of $E(K)$ representing
$\gamma\in H_1(E(K);\mz)$;
\item Assume by general position that $v$ is transversal to $L$;
\item Replace $v$ by a new field $v'$ which runs parallel to $L$ in a tubular
neighbourhood $W(L)$ of $L$; note that $\alpha(v\ristr{E(K)},v')=0\in
H_1(E(K);\mz)$;
\item Modify $v'$ only within $W(L)$ to a field $w$ which runs opposite to $L$ on
$L$
and has non-positive radial component on $W(L)$.
\end{enumerate}

\noindent Our next step is to extend $w$ to a field $z$ on the whole of $M$, which
we can do simply by defining $z$ to coincide with $v$ on $U(K)$.
Since $\gamma$ is in the kernel, at the $H_1$-level, of the inclusion of $E(K)$ into $M$,
the homotopy classes of $v$ and $z$ on $M$ differ at most by a Hopf number
({\em i.e.}~they define the same Euler structure on $M$).
Therefore we can select a ball $B$ contained in the interior of $E(K)$ and
modify $z$ on $B$ to a new field $y$ such that the Hopf number of $y$ relative
to $v$ is zero. The modification on $B$ 
is also a Pontrjagin move. Note now that $w$ and $y\ristr{E(K)}$
differ by a local modification, so they define the same Euler structure on $E(K)$.
In particular $\alpha(v\ristr{E(K)},y\ristr{E(K)})=\gamma\in H_1(E(K);\mz)$.

Now by construction $y$ and $v$ are homotopic on $M$ and $K$ is transversal to $y$. 
If $(v_t)_{t\in[0,1]}$ is the
homotopy, with $v_0=y$ and $v_1=v$, we can apply Proposition~\ref{follow:homotopy} 
and find a continuous family $(g_t)_{t\in[0,1])}$ of diffeomorphisms of $M$
fixed on $\partial M$ with $g_0={\rm id}$ and 
$g_t(K)$ transversal to $v_t$ for all $t$. 
As in the proof of Proposition~\ref{group:obstr} the homology class
$$\alpha(\xi(v,K),g_t^*(\xi(v_t,g_t(K))))$$
is constantly $\gamma$ because
it is $\gamma$ at $t=0$. So it is sufficient to define $K_\gamma$ as 
$g_1(K)$ and $f$ as $g_1^{-1}$.

When $v$ is transversal to an overtwisted contact structure
$\eta$, we fix a metric such that they are actually orthogonal, and
we modify our proof as follows
(assuming the reader is familiar with the techniques of Eliashberg, see~\cite{eliash}):
\begin{enumerate}
\item Instead of modifying $v\ristr{E(K)}$ to $w$ by a Pontrjagin move,
we construct a new contact structure by
application of a Lutz twist to $\eta\ristr{E(K)}$, so that the effect (up to
homotopy) on the orthogonal vector field is the same as the 
original modification. Then we extend the structure near $K$ as obvious,
calling $z$ its normal field.
\item Instead of modifying $z$ to $y$, again we use a Lutz twist on the
normal contact structure. Moreover we make sure that $y^{\perp}$ is overtwisted
by applying (if necessary) another Lutz twist of the sort which does not
change the homotopy class, away from $K$.
\item We conclude using Eliashberg's classification of overtwisted structures, 
according to which two such structures which are homotopic as plane fields are 
automatically isotopic.
\end{enumerate}
The proof is complete.\finedim{link:family:prop}

\begin{rem}{\em A more constructive proof of the contact version of 
Proposition~\ref{link:family:prop} may be given in the framework 
of~\cite{contspin}. On the other hand, the proof we have given above
raises the following natural
question: assume $\eta_0$ and $\eta_1$ are overtwisted contact structures on $M$,
let $L_0$ and $L_1$ be links
tangent to $\eta_0$ and $\eta_1$ respectively, and
assume there exist a family $(\eta_t,L_t)_{t\in[0,1]}$ where 
$(\eta_t)$ is a homotopy of plane fields, $(L_t)$ is an isotopy, and
$L_t$ is tangent to $\eta_t$. Can this family be replaced by a similar
one where $(\eta_t)$ is an isotopy? Eliashberg's classification theorem
may be stated as ``yes, for $L_0=\emptyset$'', and a general proof
could possibly be obtained along the lines of~\cite{eliash}.
Should the answer be ``yes, for any $L_0$'', we would have a bijection
between pseudo-Legendrian links (up to weak equivalence)
and Legendrian links in
overtwisted structures (up to Legendrian isotopy).}
\end{rem}

\paragraph{Curls and winding number.}  We show in this paragraph that torsions
are sensitive to an analogue of the winding number (the invariant which allows
to distinguish framed-isotopic link projections which are not equivalent under
the second and third of Reidemeister's moves, see~\cite{trace}). This will
allow us to give another recipe, besides  Proposition~\ref{link:family:prop},
to construct knots which are distinguished by torsion. Moreover we will give an
example of knot which is not good. The proof of the next result uses the
example of Section~\ref{exa:section}, so it is deferred to
Section~\ref{proofs}.

\begin{prop}\label{wind:sensitive}
Consider a field $v$ on $M$ and a portion of $M$ on which $v$
can be identified to the vertical field in $\mr^3$. Consider knots $K_0$ and
$K_1$ which are transversal to $v$ and differ only within the chosen portion of $M$,
as shown in Fig.~\ref{wind:fig}.
\begin{figure}
\centerline{\psfig{file=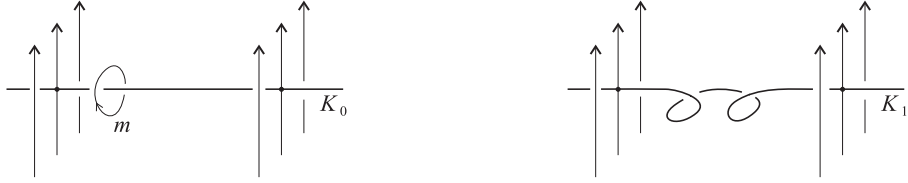,width=12cm}}
\caption{\label{wind:fig}Knots which differ for the winding number.}
\end{figure}
Choose a meridian $m$ of $K_0$ as also shown in the figure.
Let $f$ be an isotopy which maps $K_1^{(v)}$ to $K_0^{(v)}$ and is supported in
a tubular neighbourhood of $K_0$. Then:
$$\alpha(\xi(v,K_0),f_*(\xi(v,K_1)))=[m]\in H_1(E(K_0);\mz).$$
\end{prop}

\begin{prop}\label{hyper:curl} Let $(v,K_0)$ be a pseudo-Legendrian pair in
$M$, and denote by $[m]\in H_1(E(K_0);\mz)$ the homology class of the meridian of
$U(K_0)$. Assume either that $K_0^{(v)}$ is good and $[m]\neq 0$ or that
$E(K_0)$ is hyperbolic and $[m]$ has infinite order.
Let $K_1$ be a knot obtained
from $K_0$ as in Fig.~\ref{wind:fig}. Then $(v,K_0)$ and $(v,K_1)$  are not
weakly equivalent. \end{prop}

\dim{hyper:curl} By contradiction, using 
Propositions~\ref{group:obstr} and~\ref{wind:sensitive}, we would get elements
$\xi_0,\xi_1$ of $\eul(E(K_0),\pp(K_0^{(v)})$ such that $\alpha(\xi_0,\xi_1)=[m]$
and $\xi_1=f_*(\xi_0)$ for some $[f]\in G(K_0)$. If $K_0^{(v)}$ is good and $[m]\neq 0$ 
this is a contradiction. Assume now that $E(K_0)$ is hyperbolic and $[m]$ has infinite order.
Since $f_*([m])=[m]$, using 
Lemma~\ref{power:action} we easily see that $\alpha(\xi_0,f_*^k(\xi_0))=k\cdot[m]$ for
all $k$. Proposition~\ref{collar:trivial} and the result of Johansson already used in the proof of
Theorem~\ref{many:good}
now imply that $f^k$ acts trivially on $\eul(E(K_0),\pp(K_0^{(v)})$ for some
$k$, whence the contradiction.\finedim{hyper:curl}

\noindent As an application of Proposition~\ref{wind:sensitive}, 
we can show that there exist knots which are not good. Consider $S^2\times[0,1]$
with vector field parallel to the $[0,1]$ factor. Let $K_0$ be the equator
of $S^2\times\{1/2\}$, and let $K_1$ be obtained from $K_0$ by the
modification described in Fig.~\ref{wind:fig}.
Using Proposition~\ref{wind:sensitive}, if we choose a framed-isotopy
$g$ of $K_1^{(v)}$ onto $K_0^{(v)}$ supported in $U(K_0)$, we have
$$\alpha(\xi(v,K_0),(g\ristr{E(K_1)})_*(\xi(v,K_1)))=[m],$$
where $[m]$ is a generator of $H_1(E(K_0);\mz)\cong\mz$.
On the other hand, $K_1$ is strongly equivalent to $K_0$ in $(M,v)$
(the winding number only exists on $\mr^2$, not on $S^2$). 
So there exists an isotopy $h$ of $K_1^{(v)}$ onto $K_0^{(v)}$ through links
transversal to $v$, and we have
$$\alpha(\xi(v,K_0),(h\ristr{E(K_1)})_*(\xi(v,K_1)))=0.$$
This implies that $(h\compo g^{-1})\ristr{E(K_0)}$ acts non-trivially
on $\xi(v,K_0)\in\eul(E(K_0),\pp(K_0^{(v)}))$.

\section{Torsion of a convex combinatorial Euler structure}\label{torsion:def:section}
In this section we formally define torsion.
Fix a manifold $M$, a {\em convex} boundary pattern $\pp=(W,B,V,\emptyset)$ on $M$,
a cellularization $\cc$ suited to $\pp$ and
a representation $\varphi:\pi_1(M)\to\Lambda_*$, where $\Lambda$ is as mentioned 
before the statement of
Theorem~\ref{torsion:existence:statement}. We will denote by $\varphi$ again the extension
$\mz[\pi_1(M)]\to\Lambda$ (a ring homomorphism). 

We consider now the universal cover $q:\tilde{M}\to M$ and
the twisted chain complex $\Cphi_*(M,W\cup V)$, where
$\Cphi_i(M,W\cup V)$ is defined as 
$\Lambda\otimes_{\varphi}\Ccell_i(\tilde{M},q^{-1}(W\cup V);\mz)$,
and the boundary operator is induced from the ordinary boundary.
The homology of this complex is denoted by $H^\varphi_*(M,W\cup V)$ and
called the $\varphi$-twisted homology. We assume that each $H^\varphi_i(M,W\cup V)$
is a free $\Lambda$-module and fix a basis $\hbasis_i$.

\begin{rem}\label{twisted:remarks}{\em
\begin{enumerate}
\item\label{basepoint:remark} To have a formal completely intrinsic definition of  
$H^\varphi_*(M,W\cup V)$, one should fix from the
beginning a basepoint $x_0\in M$ for $\pi_1(M)$, and consider pointed
universal covers $q:(\tilde M,\tilde x_0)\to(M,x_0)$, because any two such covers
are {\em canonically} isomorphic, and the action of $\pi_1(M)$ on $\tilde{M}$ is
{\em canonically} defined on them. 
\item To define $H^\varphi_*(M,W\cup V)$
we have used in an essential way the fact that $W\cup V=\overline{W}$ is 
closed, because otherwise $\Cphi_*(M,W\cup V)$ cannot be defined.
\item\label{free:basis:remark} $\Cphi_i(M,W\cup V)$ is a free $\Lambda$-module, and
each $\mz[\pi_1(M)]$-basis of $\Ccell_i(\tilde{M},q^{-1}(W\cup V);\mz)$ determines a
$\Lambda$-basis of $\Cphi_i(M,W\cup V)$.
\item\label{varphibar:defin} If we compose $\varphi$ with the projection 
$\Lambda_*\to\Kbar_1(\Lambda)$ we get a homomorphism of $\pi_1(M)$ into an
{\em Abelian} group, so we get a homomorphism
$\overline{\varphi}:H_1(M;\mz)\to\Kbar_1(\Lambda).$
\end{enumerate}}
\end{rem}

\noindent Now let $\xi\in\eulc(M,\pp)$ and choose a representative of $\xi$ 
as in point~\ref{connected:spider:point} of 
Proposition~\ref{combin:prop:statement}, namely
$$\sum_{\sigma\in\cc,\ \sigma\subset M\setminus(W\cup V)}
\index(\sigma)\cdot \beta_\sigma$$
with $\beta_\sigma(0)=x_0$ for all $\sigma$, $x_0$ being a fixed point of $M$.
We choose $\tilde x_0\in q^{-1}(x_0)$ and consider the liftings
$\tilde\beta_\sigma$ which start at $\tilde x_0$. For $\sigma\subset M\setminus(W\cup V)$
we select its preimage $\tilde\sigma$ which contains $\tilde\beta_\sigma(1)$, and
define $\gbasis(\xi)$ as the collection of all these $\tilde\sigma$.
Arranging the $i$-dimensional elements of $\gbasis(\xi)$ in any order, by 
Remark~\ref{twisted:remarks}(\ref{free:basis:remark}) we
get a $\Lambda$-basis $\gbasis_i(\xi)$ of $\Cphi_i(M,W\cup V)$. 
We consider a set ${\tilde\hbasis}_i$ of elements of $\Cphi_i(M,W\cup V)$
which project to the fixed basis $\hbasis_i$ of $H^\varphi_i(M,W\cup V)$.

Now note that, given a free $\Lambda$-module $L$ 
and two finite bases $\basis=(b_k)$,
$\basis'=(b'_k)$ of $M$, the assumption made on $\Lambda$ guarantees that
$\basis $ and $\basis'$ have the same number of elements, so there
exists an invertible square matrix $(\lambda^h_k)$ such that
$b'_k=\sum_h\lambda^h_k b_h$. We will denote by
$[\basis'/\basis ]$ the image of $(\lambda^h_k)$ in $K_1(\Lambda)$
(see Section~\ref{statements} for the definition).

\begin{prop}\label{defining:proposition}
If $\basis_i\subset\Cphi_i(M,W\cup V)$ is such that $\partial\basis_i$ is a
$\Lambda$-basis of $\partial(\Cphi_i(M,W\cup V))$, then
$(\partial\basis_{i+1})\cdot\tilde\hbasis_i\cdot\basis_i$ is a $\Lambda$-basis of
$\Cphi_i(M,W\cup V)$, and 
$$\tau^\varphi(M,\pp,\xi,\hbasis)=\pm
\prod_{i=0}^3\Big[ \Big( (\partial\basis_{i+1})\cdot\tilde\hbasis_i\cdot
\basis_i\Big) \;\Big/\;\gbasis_i(\xi)\Big]^{(-1)^i}\in
\Kbar_1(\Lambda)$$
is independent of all choices made. Moreover 
\begin{equation}
\tau^\varphi(M,\pp,\xi',\hbasis)=\tau^\varphi(M,\pp,\xi,\hbasis)\cdot 
\overline{\varphi}(\alphac(\xi',\xi)).\label{combinatorial:equivariance:formula}
\end{equation}
\end{prop}

\dim{defining:proposition}
The first assertion and independence of the $\basis_i$'s is purely algebraic and classical,
see~\cite{milnor}. Now note that $\xi\in\eulc(M,\pp)$ was used to select the bases 
$\gbasis_i(\xi)$. The $\gbasis_i(\xi)$ are of course not uniquely determined themselves,
but we can show that different choices lead to the same value of $\tau^\varphi$.

First of all, the arbitrary ordering in the $\gbasis_i(\xi)$ is inessential because
torsion is only regarded up to sign. Second, consider the effect of choosing a
different representative of $\xi$. This leads to a new family $\tilde\sigma'$ of
cells. If $\tilde\sigma'=a(\sigma)\cdot \tilde\sigma$, with $a(\sigma)\in\pi_1(M)$,
and $\overline{a}(\sigma)$ is the image in $H_1(M;\mz)$, we automatically have
$$\sum_{\sigma\subset M\setminus{W\cup V}}\index(\sigma)\cdot \overline{a}(\sigma)=0\in
H_1(M;\mz),$$
which allows to conclude that also the representative chosen is inessential.
The choice of the lifting $\tilde x_0$ can be shown to be inessential either
in the spirit of Remark~\ref{twisted:remarks}(\ref{basepoint:remark}),
or by showing that a simultaneous $a$-translation of all $\tilde\sigma$, for
$a\in\pi_1(M)$, multiplies the torsion by $\overline{\varphi}(a)^{\chi(M)-\chi(W\cup V)}=1$.

Formula~(\ref{combinatorial:equivariance:formula}) is readily established by choosing
representatives $\sum\index(\sigma)\cdot \beta_\sigma$ and 
$\sum\index(\sigma)\cdot \beta'_\sigma$ of $\xi$ and $\xi'$ 
such that $\beta'_\sigma=\beta_\sigma$ for all 
$\sigma$ but one.\finedim{defining:proposition}

Since the above construction uses the cellularization $\cc$ in a
way which may appear to be essential,
we add a subscript $\cc$ to the torsion we have defined.
The next result, which can be established following Turaev~\cite{turaev:Euler},
shows that dependence on $\cc$ is actually inessential.
Together with Theorem~\ref{reconstruction:statement} and
Propositions~\ref{combin:prop:statement} and~\ref{defining:proposition}, it
concludes the proof of Theorem~\ref{torsion:existence:statement}.

\begin{prop}
Let $\cc$ and $\cc'$ be cellularizations suited to $\pp$. Assume that
$\cc'$ subdivides $\cc$, and consider the bijection
$\ss_{(\cc',\cc)}:\eulc(M,\pp)_\cc\to\eulc(M,\pp)_{\cc'}$ of
Proposition~\ref{combin:prop:statement}, and the canonical isomorphism
$j_{(\cc',\cc)}:H^\varphi_*(M,W\cup V)_\cc\to H^\varphi_*(M,W\cup V)_{\cc'}$. Then,
with obvious meaning of symbols we have:
$$\tau^\varphi_\cc(M,\pp,\xi,\hbasis)=
\tau^\varphi_{\cc'}(M,\pp,\ss_{(\cc',\cc)}(\xi),j_{(\cc',\cc)}(\hbasis)).$$
\end{prop}

It is maybe appropriate here to remark that the choice of a basis $\hbasis$  of
$H^\varphi_*(M,W\cup V)$ and the definition of $\tau^\varphi(M,\pp,\xi,\hbasis)$ implicitly 
assume a description of the universal cover of $M$, which is typically undoable in
practical cases. However,
if one starts from a representation of $\pi_1(M)$ into the units of a {\em commutative}
ring  $\Lambda$, {\em i.e.} a representation which factors through
 one of $H_1(M;\mz)$,
one can use from the very beginning the maximal Abelian rather
than the universal cover, which makes computations more feasible.

\begin{rem}{\em Turaev~\cite{turaev:Reidemeister} has shown that a homological
orientation yields a sign-refinement of torsion, {\em i.e.}~a lifting from 
$\Kbar_1(\Lambda)$ to $K_1(\Lambda)$. This refinement extends with minor
modifications to our setting of boundary tangency. This sign-refinement,
in the closed and monochromatic case, is an essential component of the 
theory (for instance, it is crucial for the relation with
the 3-dimensional Seiberg-Witten invariants~\cite{turaev:spinc},~\cite{turaev:nuovo}
and for the definition of the Casson invariant~\cite{lescop}), so we expect
it to be relevant also in the boundary pattern case.}
\end{rem}

\paragraph{Computation of torsion via disconnected spiders.} In this paragraph we show that 
to determine the family of lifted cells necessary to define torsion one can use representatives
of Euler structures more general than those used above. This is a technical point which we
will use below to compute torsions using branched spines (Section~\ref{spines:section}).

We fix $M$, $\pp$, $\cc$ and $\varphi$ as above, and $\xi\in\eulc(M,\pp)$. Let
$\gbasis(\xi)=\{\tilde\sigma\}$ be the family of liftings of the cells lying in
$M\setminus(W\cup V)$ determined by a connected spider as explained above. Note
that if $\gbasis'=\{\tilde\sigma'\}$ is any other family of liftings we have 
$\tilde\sigma'=a(\sigma)\cdot\tilde\sigma$ for some $a\in\pi_1(M)$, and we can
define $$h(\gbasis',\gbasis(\xi))=\sum_{\sigma\subset M\setminus (W\cup
V)}\index(\sigma)\cdot \overline{a}(\sigma)\in H_1(M;\mz).$$

\begin{prop}\label{no:need:to:lift}
Assume there exists a partition $\cc_1\sqcup\dots\sqcup\cc_k$ of the set of cells 
lying in $M\setminus(W\cup V)$, and let $\xi\in\eulc(M,\pp)$ have a representative of the form
$$z=\sum_{j=1}^k\left(\sum_{\sigma\in\cc_j\setminus\{\sigma_j\}}
\index(\sigma)\cdot \gamma^{(j)}_\sigma\right)$$
where $\sigma_j\in\cc_j$ and $\gamma^{(j)}_\sigma:([0,1],0,1)\to(M,p_{\sigma_j},p_\sigma)$.
Choose any lifting $\tilde p_{\sigma_j}$ of $p_{\sigma_j}$, lift 
$\gamma^{(j)}_\sigma$ to $\tilde\gamma^{(j)}_\sigma$ starting from $\tilde p_{\sigma_j}$,
let $\tilde\sigma'$ be the lifting of $\sigma$ containing $\tilde\gamma^{(j)}_\sigma(1)$,
and let $\gbasis'$ be the family of all these liftings.
Then $h(\gbasis',\gbasis(\xi))=0\in H_1(M;\mz)$. In particular $\gbasis'$ can be used to compute 
$\tau^\varphi(M,\pp,\xi,\hbasis)$.
\end{prop}

\dim{no:need:to:lift}
Note first that the coefficient of $p_{\sigma_j}$ in $\partial z$ is exactly
$$-\sum_{\sigma\in\cc_j\setminus\{\sigma_j\}}\index(\sigma).$$
On the other hand this coefficient must be equal to $\index(\sigma_j)$. 
Summing up we deduce that $\sum_{\sigma\in\cc_j}\index(\sigma)=0$.

Now choose $x_0\in M$ and $\delta^{(j)}:([0,1],0,1)\to(M,x_0,p_{\sigma_j})$. For 
$\sigma\in\cc_j$ define
$$\beta_\sigma=\cases{
\delta^{(j)} & if $\sigma=\sigma_j$ \cr
\delta^{(j)}\cdot\gamma^{(j)}_\sigma & otherwise,}$$
so that $\beta_\sigma:([0,1],0,1)\to(M,x_0,p_\sigma)$, whence 
$w=\sum_{\sigma\subset M\setminus(W\cup V)}\beta_\sigma$ is an Euler chain.
Moreover:
$$w-z=\sum_{j=1}^k\left(\sum_{\sigma\in\cc_j}\index(\sigma)\right)
\cdot\delta^{(j)}=0\in H_1(M;\mz),$$
so $[w]=\xi$. Now choose $\tilde x_0$ over $x_0$, lift the $\delta^{(j)}$ and 
$\beta_\sigma$ starting from $\tilde x_0$, and let
$a^{(j)}\in\pi_1(M)$ be such that $\tilde p_{\sigma_j}=a^{(j)}\cdot\tilde\delta^{(j)}(1)$.
Then
$$h(\gbasis',\gbasis(\xi))=\sum_{j=1}^k\left(\sum_{\sigma\in\cc_j}\index(\sigma)\right)
\cdot\overline{a}^{(j)}=0\in H_1(M;\mz),$$
and the proof is complete.\finedim{no:need:to:lift}

\section{Spines and computation of torsion}\label{spines:section}
In this section we show how to compute torsions starting from a combinatorial
encoding of vector fields. We first review
the theory developed in~\cite{lnm}. See the beginning of Section~\ref{statements}
for our conventions on manifolds, maps, and fields. In addition to the terminology
introduced there, we will need the notion of {\em traversing} field on a manifold
$M$, defined as a field whose orbits eventually intersect $\partial M$ 
transversely in both directions (in other words, orbits are compact intervals).

\paragraph{Branched spines.}
A {\em simple} polyhedron $P$ is a  finite connected 2-dimensional polyhedron
with singularity of stable nature (triple lines and points where six non-singular
components meet). Such a $P$ is called {\it standard} if all the components of
the natural stratification given by singularity are open cells. Depending on
dimension, we will call the components {\it vertices, edges} and {\it regions}.

\begin{figure}
\centerline{\psfig{file=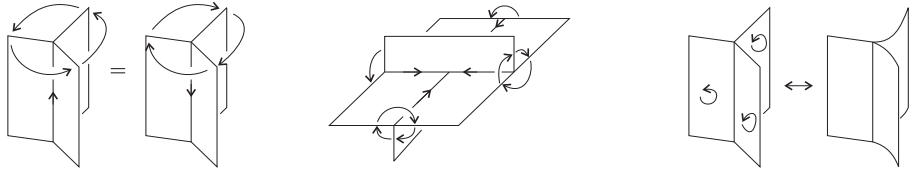,width=12cm}}
\caption{\label{screw:branch} Convention on screw-orientations, compatibility 
at vertices, and geometric interpretation of branching.}
\end{figure}

A {\em standard spine} of a $3$-manifold $M$ with $\partial M\neq\emptyset$
is a standard polyhedron $P$ embedded in $\interior(M)$ so that $M$ collapses onto $P$.
Standard spines of oriented $3$-manifolds are characterized among standard polyhedra
by the property of carrying an {\em orientation}, defined 
(see Definition~2.1.1 in~\cite{lnm}) as a ``screw-orientation''
along the edges (as in the left-hand-side of Fig.~\ref{screw:branch}),
with an obvious compatibility at vertices
(as in the centre of Fig.~\ref{screw:branch}).
It is the starting point of the theory of standard spines
that every oriented $3$-manifold $M$ with $\partial M\neq\emptyset$ has an oriented
standard spine, and can be reconstructed (uniquely up to homeomorphism) 
from any of its oriented standard spines. See~\cite{casler} for the non-oriented version 
of this result and~\cite{manuscripta} or Proposition~2.1.2 in~\cite{lnm} for the (slight) oriented refinement.

A {\it branching} on a standard polyhedron $P$ is an
orientation for each region of $P$, such that no edge is induced the same
orientation three times. See the right-hand side of Fig.~\ref{screw:branch}
and Definition 3.1.1 in~\cite{lnm} for the geometric meaning of this notion.
An oriented standard spine $P$ endowed with a branching is shortly named 
{\em branched spine}. We will never use specific notations for the extra structures: 
they will be considered to be part of $P$.
The following result, proved as Theorem~4.1.9 in~\cite{lnm}, is the starting point of our constructions.

\begin{prop}\label{from:spine:to:field}
To every branched spine $P$ there corresponds a manifold $M(P)$ 
with non-empty boundary and a concave traversing field $v(P)$ on $M(P)$.
The pair $(M(P),v(P))$ is well-defined up to diffeomorphism.
Moreover an embedding $i:P\to\interior(M(P))$ is defined, 
and has the property that $v(P)$ is positively transversal to $i(P)$.
\end{prop}

\begin{figure}
\centerline{\psfig{file=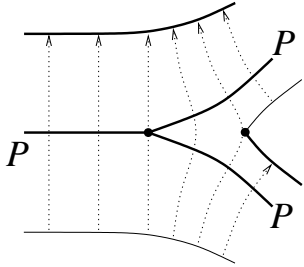,width=4cm}}
\caption{\label{constr:M} Manifold and field associated to a branched spine.}
\end{figure}

The topological construction which underlies this proposition is actually quite
simple, and it is illustrated in Fig.~\ref{constr:M}. Concerning the last
assertion of the proposition, note that the branching allows to define an
oriented tangent plane at each point of $P$.

\paragraph{Combinatorial encoding of combings.}
Let $P$ be a branched spine, and define $v(P)$ on $M(P)$ as just explained.
Assume that in $\partial M(P)$ there is only one component which is homeomorphic to $S^2$
and is split by the tangency line of $v(P)$ to $\partial M(P)$ into two
discs. (Such a component will be denoted by $\stwotriv$.) Now, notice that $\stwotriv$ is
also the boundary of the closed $3$-ball with constant vertical field, 
denoted by $B^3_{\rm triv}$. This shows that we can cap off
$\stwotriv$ by attaching a copy of $B^3_{\rm triv}$, 
getting a compact manifold $\hatM(P)$ and a field $\hatv(P)$ on $\hatM(P)$. If we denote by $\hatpp(P)$ the boundary pattern of $\hatv(P)$ on
$\hatM(P)$, we easily see that the pair $(\hatM(P),\hatv(P))$ is only well-defined
up to homeomorphism of $\hatM(P)$ and homotopy of $\hatv(P)$ through fields compatible
with $\hatpp(P)$. Note also that $\hatpp(P)$ is automatically concave.

If $\pp$ is a boundary pattern on $M$, we define ${\rm Comb}(M,\pp)$ as
the set of fields compatible with $\pp$ under homotopy through fields also
compatible with $\pp$. An element of ${\rm Comb}(M,\pp)$
is called a {\em combing} on $(M,\pp)$. Note that we have a projection
${\rm Comb}(M,\pp)\to\eul(M,\pp)$.

The above construction shows that
a branched spine $P$ with only one $\stwotriv$ on $\partial M(P)$ defines
an element $\Phi(P)$ of ${\rm Comb}(\hatM(P),\hatpp(P))$.
One of the main achievements of~\cite{lnm} (Theorems~1.4.1 and~5.2.1)
is the following.

\begin{teo}\label{closed:bijectivity:teo}
\begin{enumerate}
\item\label{closed:surg} If $M$ is a closed oriented $3$-manifold, $\Phi$ maps surjectively
the set $\{P:\hatM(P)\cong M\}$ onto ${\rm Comb}(M)={\rm Comb}(M,\emptyset)$.
\item\label{closed:inj} A finite list of local combinatorial moves on branched spines
can be given so that if $\hatM(P_0)\cong\hatM(P_1)\cong M$ is closed and 
$\Phi(P_0)=\Phi(P_1)\in {\rm Comb}(M)$, then $P_1$ is obtained from $P_0$ by a
finite sequence of these moves.
\end{enumerate}
\end{teo}

\begin{figure}
\centerline{\psfig{file=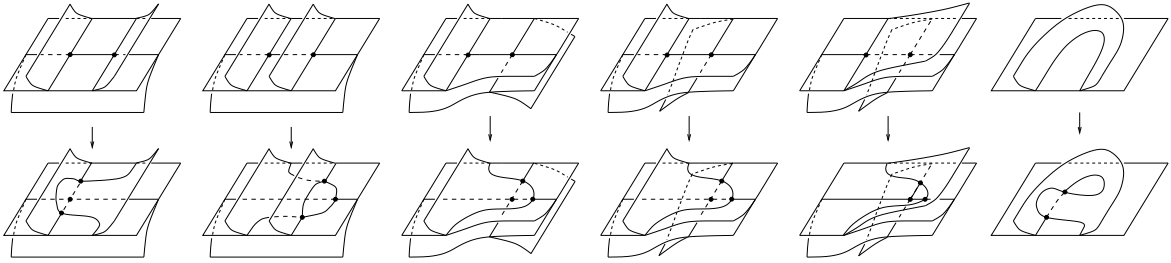,width=15.5cm}}
\caption{\label{phys:moves} Moves for branched standard spines.}
\end{figure}

\noindent In the present paper we will not use the moves referred to in the
previous statement, but to give the reader an idea of their geometric meaning
we quickly picture them.
The complete list actually consists of 18 moves, but the essential ``physical''
phenomena which occur are only those shown in 
Fig.~\ref{phys:moves} (the other moves are obtained by taking mirrors of those shown).

In~\cite{second:paper} we will show that the rightmost move in 
Fig.~\ref{phys:moves} is actually implied by the other moves, and we will
establish the following extension of Theorem~\ref{closed:bijectivity:teo}.

\begin{teo}\label{bounded:bijectivity:teo}
\begin{enumerate}
\item\label{bounded:surg} If $M$ is any compact oriented 
$3$-manifold and $\pp$ is a concave boundary pattern
on $M$ not containing $\stwotriv$ components, then $\Phi$ maps surjectively
$\{P:\ \hatM(P)\cong M,\ \hatpp(P)\cong\pp\}$ onto ${\rm Comb}(M,\pp)$.
\item\label{bounded:inj} The same finite list of moves as in point~\ref{bounded:inj} of
Theorem~\ref{closed:bijectivity:teo} has the property that if 
$(\hatM(P_0),\hatpp(P_0))\cong(\hatM(P_1),\hatpp(P_1))\cong(M,\pp)$ is as above and 
$\Phi(P_0)=\Phi(P_1)\in {\rm Comb}(M,\pp)$, then $P_1$ is obtained from $P_0$ by a
finite sequence of these moves.
\end{enumerate}
\end{teo}

\noindent The proof of this result requires considerable technicalities, so we have
decided to omit it here, also because point~\ref{bounded:inj} is not used, 
and point~\ref{bounded:surg} is only needed to show that the recipe we will
give to compute torsions actually allows to compute {\em all} concave torsions.
We just mention that both points are established by extending to the bounded
case the notion of {\em normal section} of a field, introduced in~\cite{ishii}
and~\cite{lnm} (Section 5.1). The following
geometric interpretation of point~\ref{bounded:surg} may be of some interest.

\begin{rem}\label{trivial:pieces}
{\em In general, the dynamics of a field, even a concave one, can be very complicated,
whereas the dynamics of a traversing field (in particular, $B^3_{\rm triv}$) is simple.
Point~\ref{bounded:surg} in Theorem~\ref{bounded:bijectivity:teo} means that for any 
(complicated) concave field
there exists a sphere $S^2$ which splits the field into two (simple) pieces:
a standard $B^3_{\rm triv}$ and a concave traversing field.}
\end{rem}

Another reason for not proving point~\ref{bounded:surg} of
Theorem~\ref{bounded:bijectivity:teo} in general is that we can give an easy special 
proof for the case we are most interested in, namely link complements.
Note that our argument relies on 
Theorem~\ref{closed:bijectivity:teo} (and its proof).

\vspace{1pt}
\noindent{\em Proof of point~\ref{bounded:surg} of Theorem~\ref{bounded:bijectivity:teo}
for link complements}. We have to show that if $M$ is closed, $v$ is a field on $M$
and $L$ is transversal to $v$, then the complement $E(L)$ of $L$ with the restricted field
is represented by some branched spine in the sense explained above.

The construction explained in Section~5.1  of~\cite{lnm} shows that 
there exists a branched standard spine $P$ such that $v$ is positively 
transversal to $P$ and the complement of $P$, with the 
restriction of $v$, is isomorphic to
the open 3-ball with the constant vertical field. 
The last condition easily implies that $L$ can be 
isotoped through links transversal to $v$ to a link lying in an arbitrarily small 
neighbourhood of $P$, with the further property that its natural projection
on $P$ is $\cont^1$, possibly with crossings. This fact is the starting 
point of a treatment of framed links via $\cont^1$ projections
on spines, which we plan to develop in~\cite{second:paper}.

\begin{figure}
\centerline{\psfig{file=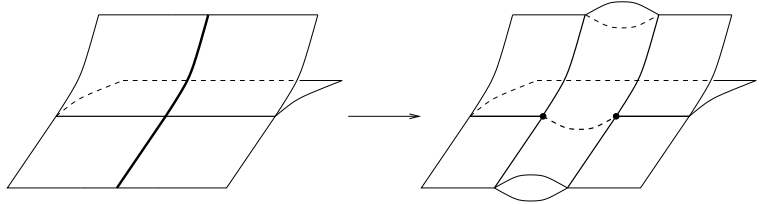,width=10cm}}
\caption{\label{dig:tunnel} How to dig a tunnel in a spine.}
\end{figure}

Once $L$ has been isotoped to a $\cont^1$ link on $P$,
a branched spine of $(E(L),v\ristr{E(L)})$ is obtained
by digging a tunnel in $P$ along the projection of $L$, as shown in Fig.~\ref{dig:tunnel}.
A crossing in the projection will of course give rise to 4 vertices in the spine.
Note that the spine which results from the digging may occasionally be 
non-standard, but it is standard as soon as the projection is complicated enough 
({\em e.g.}~if on each component there are both a crossing and an intersection with $S(P)$).
{{\hfill\hbox{\enspace\fbox{\ref{bounded:bijectivity:teo}(\ref{bounded:surg})$E(L)$}}}\vspace{5pt}}

\begin{rem}\label{no:crossing:label}
{\em Using the fact that all the regions of a branched spine
$P$ have non-empty boundary one can
show quite easily that a link $L$ with 
$\cont^1$ projection on $P$ can be isotoped through links transversal to $v(P)$ to
a link whose projection does not have crossings. An example of how to get
\begin{figure}
\centerline{\psfig{file=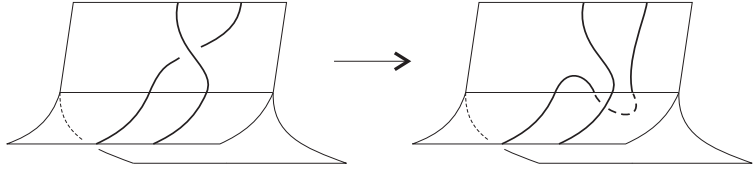,width=10cm}}
\caption{\label{rid:cross} Removing a crossing from a $\cont^1$-projection.}
\end{figure}
rid of a crossing is given in Fig.~\ref{rid:cross}.}
\end{rem}

\paragraph{Spines and ideal triangulations.} We remind the reader that an {\it
ideal triangulation} of a manifold $M$ with non-empty boundary is  a partition
$\calt$ of $\interior(M)$ into open cells of dimensions 1, 2 and 3, induced by a
triangulation $\calt'$ of the space $Q(M)$, where:
\begin{enumerate}
\item $Q(M)$ is obtained from $M$ by collapsing each component of
$\partial M$ to a point;
\item $\calt'$ is a triangulation only in a loose sense, namely
self-adjacencies and multiple adjacencies of tetrahedra are allowed;
\item The vertices of $\calt'$ are precisely the points of $Q(M)$
which correspond to components of $\partial M$.
\end{enumerate}

\begin{figure}
\centerline{\psfig{file=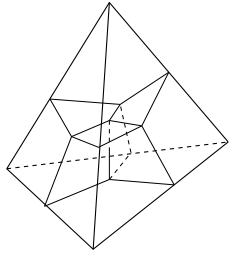,width=3cm}}
\caption{\label{duality} Duality between standard spines and ideal triangulations.}
\end{figure}

\noindent It turns out (see for instance~\cite{mafo},~\cite{tesi},~\cite{matv:new}) that there exists a natural
bijection between standard spines and ideal triangulations
of a 3-manifold. Given an ideal triangulation, the
corresponding standard spine is just the 2-skeleton of the dual
cellularization, as illustrated in Figure~\ref{duality}.
The inverse of this correspondence will be denoted by $P\mapsto\calt(P)$. 

Now let $P$ be a branched spine. First of all, we can realize $\calt(P)$ in such
a way that its edges are orbits of the restriction of $v(P)$ to
$\interior(M(P))$, and the 2-faces are unions of such orbits.  Being orbits,
the edges of $\calt(P)$ have a natural orientation, and the branching condition,
as remarked in~\cite{gr:1},
is equivalent to the fact that on each tetrahedron of $\calt(P)$ 
exactly one of the vertices is a sink and one is a source. 

\begin{rem}\label{all:oriented}
{\em It turns out that if $P$ is a branched spine, 
not only the edges, but also the faces
and the tetrahedra of $\calt(P)$ have natural orientations. For tetrahedra, we
just restrict the orientation of $M(P)$. For faces, we first note that the edges
of $P$ have a natural orientation (the prevailing orientation induced by the
incident regions). Now, we orient a face of $\calt(P)$ so that the algebraic
intersection in $M(P)$ with the dual edge is positive. }
\end{rem}

\paragraph{Euler chain defined by a branched spine.}
We fix in this paragraph a standard spine $P$ and consider its manifold $M=M(P)$.
We start by noting that the ideal triangulation $\calt=\calt(P)$ defined by $P$ 
can be interpreted as a realization of $\interior(M)$ by face-pairings on a finite
set of tetrahedra with vertices removed. If, instead of 
removing vertices, we remove open conic neighbourhoods of the vertices, thus getting {\em truncated} tetrahedra, after the 
face-pairings we obtain $M$ itself. This shows that $P$ determines a cellularization
$\caltbar=\caltbar(P)$ 
of $M$ with vertices only on $\partial M$ and 2-faces which are either triangles
contained in $\partial M$ or hexagons contained in $\interior(M)$, with edges contained alternatively in $\partial M$ and in $\interior(M)$.

\begin{figure}
\centerline{\psfig{file=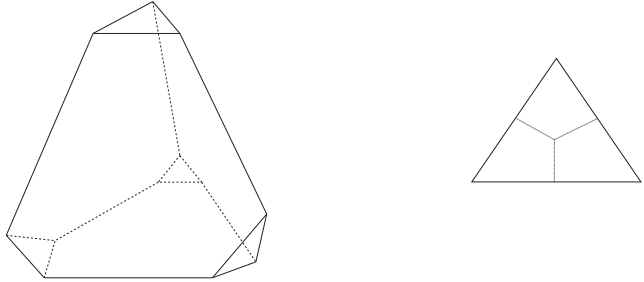,width=8.5cm}}
\caption{\label{trunc:tetra}Truncated tetrahedra and subdivision of the triangles on the
boundary}
\end{figure}

Now assume that $P$ is branched and 
that $\partial M$ contains only one $\stwotriv$ component, so 
$\hatM=\hatM(P)$ is defined. Note that $\hatM$ can be thought of as the space obtained from
$M$ by contracting $\stwotriv$ to a point, so a projection $\pi:M\to\hatM$ is defined,
and $\pi(\caltbar)$ is a cellularization of $\hatM$. Next, we modify 
$\pi(\caltbar)$ by subdividing the triangles on
$\partial\hatM$ as shown in Fig.~\ref{trunc:tetra}.
The result is a cellularization $\hattt=\hattt(P)$ of $\hatM$. 
Note that $\hattt$ on $\partial\hatM$ 
consists of ``kites'', with long edges coming from tetrahedra and short edges coming from 
subdivision. Note also that $\hattt$ has exactly one vertex $x_0$ in $\interior(\hatM)$, and
that the cells contained in $\interior(\hatM)$, except $x_0$, are the duals to the 
cells of the natural cellularization $\uu=\uu(P)$ of $P$. For $u\in\uu$ we denote by $\hat{u}$ 
its dual and by $p_u=p_{\hat{u}}$ the point where $u$ and $\hat{u}$ intersect, called the {\em centre} of both.

We will now use the field $\hatv=\hatv(P)$ to construct a combinatorial Euler chain
on $\hatM$ with respect to
$\hattt$. It is actually convenient to consider, instead of $\hatv$, the field 
$\vbar=\pi(v)$, which 
coincides with $\hatv$ except near $x_0$, where it has a 
(removable) singularity. For $u\in\uu$ we denote by $\beta_u$ the arc obtained by integrating $\vbar(P)$ in the positive direction, starting from $p_u$, until the 
boundary or the singularity is reached. We define:
$$s(P)=\sum_{u\in\uu}\index(u)\cdot \beta_u.$$

\begin{figure}
\centerline{\psfig{file=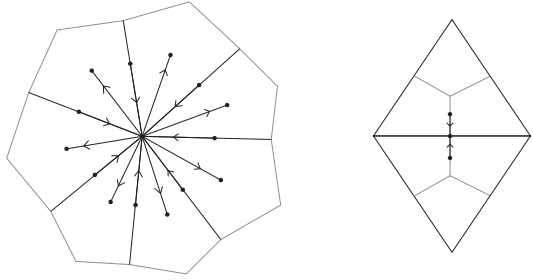,width=7cm}}
\caption{\label{starfig}The star ${\rm St}(p)$ centred at a vertex $p$ contained in $B$
and the bi-arrow ${\rm Ba}(\sigma)$ based at the midpoint of an edge $\sigma$ contained in $B$}
\end{figure}

Let us consider now the pattern $\hatpp=\hatpp(P)=(W,B,\emptyset,C)$ defined by $P$.
If $p$ is a vertex of $\pi(\caltbar)$ contained in $B$, we define its star ${\rm St}(p)$ as
the sum of the straight segments going from $p$ to the centres of all the kites containing $p$,
minus the sum of the straight
segments going from $p$ to the centres of all the long edges 
containing $p$. If $\sigma$ is an edge of $\pi(\caltbar)$ contained in $B$ we define its
bi-arrow ${\rm Ba}(\sigma)$ as the sum of the two straight segments going from the centre
$p_\sigma$ of $\sigma$ to the centres of the two short kite-edges containing $p_\sigma$.
A star and a bi-arrow are shown in Fig.~\ref{starfig}.
We define:
$$s'(P)=s(P)+\sum_{p\in B\cap\caltbar(P)^{(0)}}{\rm St}(p)+
\sum_{\sigma\in\hattt(P)^{(1)},\sigma\subset B}{\rm Ba}(\sigma).$$

\begin{lem}\label{s:prime:Euler}
$s'(P)$ defines an element of $\eulc(\hatM,\theta(P))$.
\end{lem}

\dim{s:prime:Euler}
Recall that $\theta(\hatpp)=(W,B,C,\emptyset)$, {\em i.e.}~the concave line $C$ is turned into
a convex one. So by definition we have to show that $\partial s'(P)$ contains, with the right 
sign, the centres of all cells of $\hattt$ except those of $W\cup C$. 

It will be convenient to analyze first the natural lifting of $s(P)$ to $M$, denoted by
$\tilde s(P)=\sum_{u\in\uu}\index(u)\cdot \tilde\beta_u$ with obvious meaning of symbols. So 
\begin{equation}
\partial\tilde s(P)=\sum_{u\in\uu}-\index(u)\cdot \tilde\beta_u(0)+
\sum_{u\in\uu}\index(u)\cdot \tilde\beta_u(1).
\label{s:tilde:boundary}
\end{equation}

Since the cellularization $\caltbar$ of $M$ is dual to $\uu$, the first half of 
(\ref{s:tilde:boundary}) gives the centres of the cells contained in $\interior(M)$, with right 
sign. One easily sees that the second half gives exactly the centres of the cells 
(of $\caltbar$) contained in $B$, also with right sign. 

When we project to $\hatM$ and consider $\partial s(P)$, the first half of 
(\ref{s:tilde:boundary}) again provides (with right sign)
the centres of the all cells contained in $\interior(\hatM)$, except the special
vertex $x_0$ obtained by collapsing $\stwotriv$. We can further split
the points of the second half of (\ref{s:tilde:boundary}) into those which lie on $\stwotriv$
and those which do not. The points of the first type project to $x_0$, and the resulting coefficient of $x_0$ is $\chi(B\cap\stwotriv)$, but $B\cap\stwotriv$ is an open 2-disc,
so the coefficient is 1. (We are here using the very special property of dimension 2
that $\chi$ can be computed using a finite cellularization of an open manifold, because the
boundary of the closure has $\chi=0$.) The points of the second type faithfully project to 
$\hatM$, giving the centres of the simplices contained in $B$ of the triangulation
$\pi(\caltbar)\ristr{\partial\hatM}$.
However $\hattt$ on $\partial\hatM$ is a subdivision of $\pi(\caltbar)$,
and this is the reason why we have added the stars and the bi-arrows to $s(P)$ getting $s'(P)$.
The following computation of the coefficients in $\partial s'(P)$ of the centres of the 
cells of $\hattt$ contained in $B$ concludes the proof.
\begin{enumerate}
\addtocounter{enumi}{-1}
\item Cells of dimension 0 are listed as follows:
\begin{enumerate}
\item Centres of triangles of $\pi(\caltbar)$, 
which receive coefficient $+1$ from $\partial s(P)$;
\item Midpoints of edges of $\pi(\caltbar)$, 
which receive coefficient $-1$ from $\partial s(P)$
and $+2$ from the bi-arrows they determine;
\item Vertices of $\pi(\caltbar)$ receive $+1$ from $\partial s(P)$ and
(algebraically) $0$  from the star they determine;
\end{enumerate}
\item Cells of dimension 1 are:
\begin{enumerate}
\item Short edges of kites, whose midpoints receive $-1$ from the bi-arrows;
\item Long edges of kites, whose midpoints receive $-1$ from the stars;
\end{enumerate}
\item Cells of dimension 2 are kites, and their centres receive $+1$ from the stars.
\end{enumerate}
\finedim{s:prime:Euler}

\noindent Now we denote by $\gamma_j:(0,1)\to C$, for $j=1,\dots,n$, 
orientation-preserving parameterizations of the 1-cells of $\hattt$ contained in $C$, and we extend the $\gamma_j$ to $[0,1]$, without changing notation. We define
$$s''(P)=s'(P)+\sum_{j=1}^n\gamma_j\ristr{[1/2,1]}.$$ 

\begin{lem}\label{s:second:Euler}
$s''(P)$ defines an element of $\eulc(\hatM,\hatpp)$, and 
$$[s'(P)]=\Thetac([s''(P)])\in \eulc(\hatM,\theta(\pp)).$$
\end{lem}

\dim{s:second:Euler} At the level of representatives, the second assertion is obvious,
and it implies the first assertion.
\finedim{s:second:Euler} 

\noindent We defer to Section~\ref{proofs} the proof of the next result, which shows that
the map $P\mapsto[s''(P)]\in\eulc(\hatM,\hatpp)$ allows,
using branched spines, to explicitly find the inverse
of the reconstruction map $\Psi$ of 
Theorem~\ref{reconstruction:statement}. This result was informally announced as
Theorem~\ref{informal:algorithmic} in the Introduction.

\begin{teo}\label{spider:structure:teo}
$\Psi([s''(P)])=[\hatv(P)]\in\euls(\hatM,\hatpp)$.
\end{teo}

Recall now that we have defined torsions directly only for convex patterns,
and we have extended the definition to concave patterns via the map $\Theta$.
As a consequence of Lemma~\ref{s:second:Euler} and Theorem~\ref{spider:structure:teo},
and by direct inspection of $s'(P)$, we have the following result which summarizes
our investigations on the relation between spines and torsion:

\begin{teo}\label{spines:compute:teo}
If $P$ is a branched spine which represents a manifold $\hatM$ with 
concave boundary pattern $\hatpp=(W,B,\emptyset,C)$
in the sense of Theorem~\ref{bounded:bijectivity:teo}(\ref{bounded:surg}), then for any 
representation $\varphi:\pi_1(M)\to\Lambda_*$ and any $\Lambda$-basis
$\hbasis$ of $H^\varphi_*(\hatM,W\cup C)$, the torsion
$\tau^\varphi(\hatM,\hatpp,[\hatv(M)],\hbasis)$ can be computed using
(in the sense of Proposition~\ref{no:need:to:lift}) the lifting
to the universal cover of $\hatM$ of the chain $s'(P)$ defined above. In particular, 
$s'(P)$ can be used directly, without replacing it by a connected spider.
\end{teo}

\paragraph{Boundary operators.}
To actually compute torsion starting from a branched spine $P$, besides describing
the universal (or maximal Abelian) cover of $\hatM=\hatM(P)$ and determining the preferred
liftings of the cells in $\hatM\setminus(W\cup C)$, one needs to compute
the boundary operators in the twisted chain complex 
$\Cphi_*(M,W\cup C)$. These operators are of course twisted liftings of the
corresponding operators in the cellular chain complex of $(\hatM,W\cup C)$,
with respect to $\hattt$. We briefly describe here the form of the latter operators.
Recall first that $\hattt$ consists of a special vertex $x_0$, the kites
(with their vertices and edges) on $\hatM$, and the duals of the
cells of $P$. On $\partial\hatM$ the situation is easily described, so we consider the
internal cells. 
\begin{enumerate}
\item If $R$ is a region of $P$, the ends of its dual edge 
$\hat{R}$ are either $x_0$ or vertices of $\partial\hatM$ 
contained only in long edges of kites. 
\item If $e$ is an edge of $P$ then
$\partial\hat{e}$ is given by
$\hat{R}_1+\hat{R}_2-\hat{R}_0$ plus 3 long edges of kites, where $R_0,R_1,R_2$ are the
regions incident to $e$, numbered so that $R_1$ and $R_2$ induce on $e$ the same
orientation. Here $R_0,R_1,R_2$ need not be different from each other, so
the formula may actually have some cancelation. The 3 long edges
of kites must be given an appropriate sign, and some of them 
may actually be collapsed to the point $x_0$. Note that we have only 3 kite-edges, out
of the 6 which geometrically appear on $\partial\hat{e}$, because the other 3 are white.
\item If $v$ is a vertex of $P$ then
$\partial\hat{v}$ is given by 
$\hat{e}_1+\hat{e}_2-\hat{e}_3-\hat{e}_4$ plus 6 kites, where $e_1,e_2$ are
the edges which (with respect to the natural orientation) are leaving $v$, and
$e_3,e_4$ are those which are reaching it. Again, there could be
repetitions in the $e_i$'s. The kites all have coefficient $+1$, and again some of them may
actually be collapsed to $x_0$. As above, we have only
6 kites because the other 6 are white.
\end{enumerate}

\begin{rem}\label{do:not:cut}{\em To define the cellularization $\hattt(P)$ associated to 
a spine we have decided to subdivide all the triangles on $\partial\hatM$
into 3 kites, but when doing actual computations this is not necessary and 
impractical. The only triangles 
which we really need to subdivide are those intersected by $C$, because we need the 
cellularization to be suited to the pattern. Let us consider the 4
triangles corresponding to the ends of a certain tetrahedron. If in each of them
we count the number of black kites and the number of white kites, we get respectively
$(3,0)$, $(2,1)$, $(1,2)$, $(0,3)$. So, the first and last triangles do not have to
be subdivided, and the other two can be subdivided using a segment only. Summing up,
for each vertex of $P$ we only need to add two segments on the boundary. Before
projecting $M(P)$ to $\hatM(P)$ one sees that the number of cells, with respect to
$\caltbar(P)$, is increased in all dimensions 0, 1 and 2 by twice the number of vertices of $P$.
When projecting to $\hatM(P)$ the cells lying in $\stwotriv$ get collapsed to points.}
\end{rem}

\section{An example}\label{exa:section}
\begin{figure}
\centerline{\psfig{file=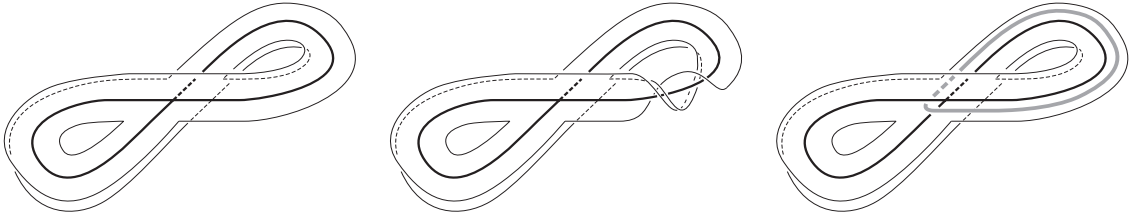,width=15.5cm}}
\caption{\label{abalone}The abalone, and a $\cont^1$ knot on it.}
\end{figure}

Figure~\ref{abalone} shows a neighbourhood of the singular set of the so-called
abalone, a branched standard spine of $S^3$, which we denote by $A$. Note that $A$ has one vertex,
two edges and two regions. The figure on the left is easier to understand, but 
it does not represent the genuine embedding of $A$ in $S^3$, which is instead
shown in the centre (hint: compute linking numbers). On the right we show (using the
easier picture) a $\cont^1$ knot $K$ on $A$. Using the genuine picture one sees that
$K$ is actually trivial in $S^3$, and its framing is $-1$.
So the knot complement $E(K)$ is actually a solid torus,
with an induced Euler structure $\xi$, and the white annulus $W\subset\partial E(K)$
is a longitudinal one. Let us now take the representation $\varphi:\pi_1(E(K))\to\mz[t^{\pm1}]$
which maps the generator to $t$. It is not hard to see that $H_*^\varphi(E(K),\overline{W})=0$, so we can
compute $\tau^\varphi(E(K),\xi)$. We describe the method to be followed, skipping several details
and all explicit formulae.

\begin{figure}
\centerline{\psfig{file=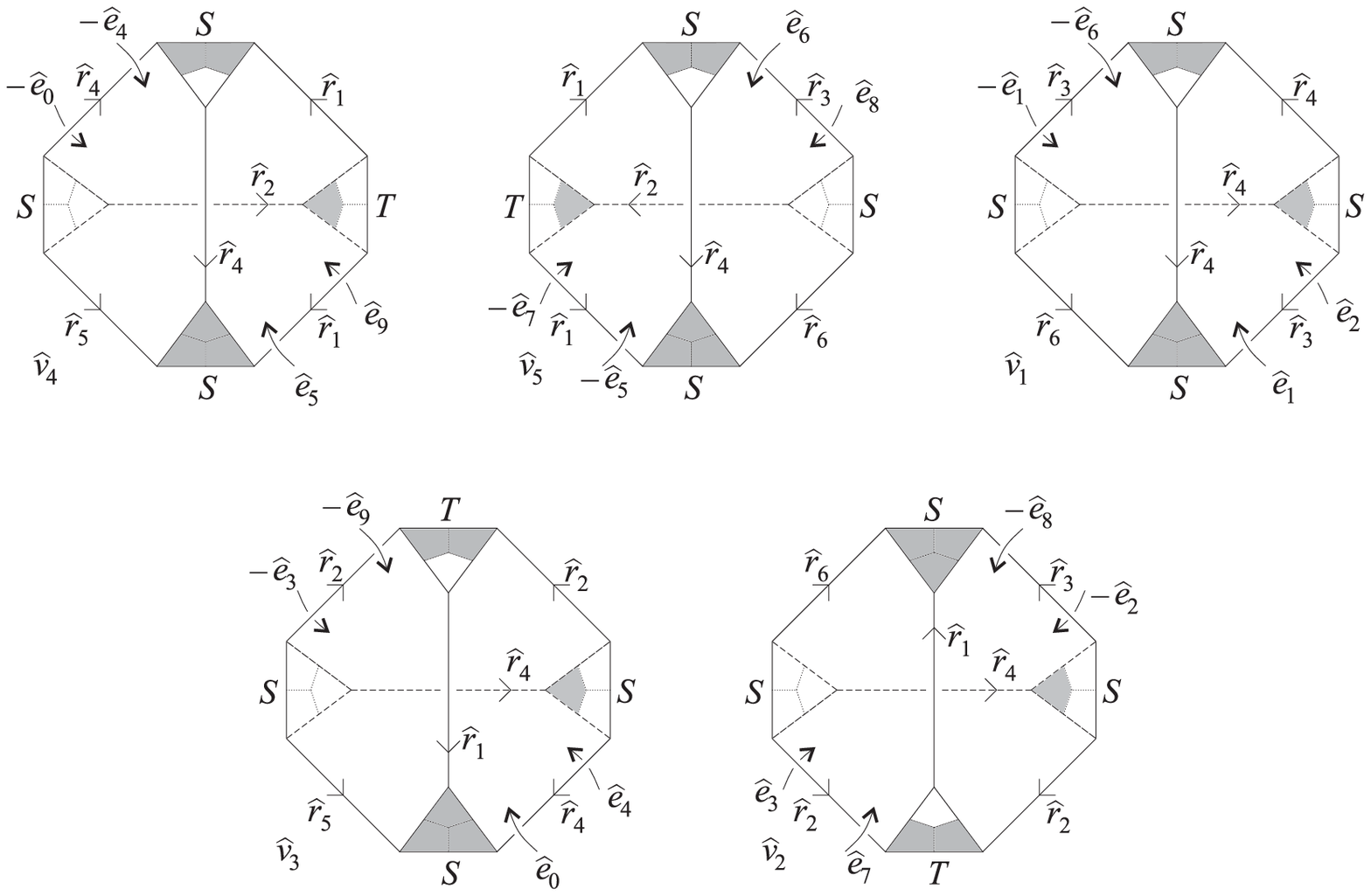,width=15.5cm}}
\caption{\label{dualtria}Truncated ideal triangulation of the knot complement.}
\end{figure}

We can apply directly the method described in the (partial) proof of 
Theorem~\ref{bounded:bijectivity:teo}, to get a branched standard spine $P$
(in the sense of Theorem~\ref{bounded:bijectivity:teo}) of $E(K)$.
This $P$ is easily recognized to have 5 vertices (denoted $v_1,\dots,v_5$),
10 edges (denoted $e_0,\dots,e_9$) and 6 regions (denoted $r_1,\dots r_6$).
Figure~\ref{dualtria} shows the truncated ideal triangulation dual to $P$.
In the figure the hat denotes duality as usual. We have written
$-\hat e_i$ instead of $\hat e_i$ when $\hat e_i$ lies on 
$\hat v_j$ but the natural
orientation of $\hat e_i$ is not induced by the orientation of $\hat v_j$.
The letters $S$ and $T$ refer to the boundary sphere and torus respectively 
($S$ should actually be collapsed to one point $x_0$, but the picture is
easier to understand before collapse).

Recall that the algebraic complex of which we must compute the
torsion has one generator for each cell in the cellularization 
of $E(K)$ arising from $P$,
excluding the white cells and the tangency circles on the boundary.
From Fig.~\ref{dualtria} we can see how many such cells there will
be in each dimension, namely 3 in dimension 0 ($x_0$ and two vertices
on $T$), 14 in dimension 1 (the $\hat r_i$'s and 8 edges on $T$),
16 in dimension 2 (the $\hat e_i$'s and the 6 black kites on $T$)
and 5 in dimension 3 (the $\hat v_i$'s). We can also easily describe
the combinatorial Euler chain $s'(P)$ which will be used to find
the preferred cell liftings: besides the orbits of the field there are
only one star and one bi-arrow; the support of $s'(P)$ has 3 connected 
components (one spider with 19 legs and head at $x_0$, the star union
the second half of $\hat r_2$, and the bi-arrow union a segment 
contained in $\hat e_3$).

To actually determine the preferred liftings we need an effective
description of the lifting of the cellularization
to the universal cover $\tilde E(K)\to E(K)$.
Since $\pi_1(E(K))=\mz$, each cell $c$ will have liftings $\tilde c^{(n)}$
for $n\in\mz$, where $\tilde c^{(n)}$ is the $n$-th translate of $\tilde c^{(0)}$.
The choice of $\tilde c^{(0)}$ itself is of course arbitrary, but to
understand the cover we must express the $\partial\tilde c^{(0)}$'s in terms of
the other $\tilde d^{(n)}$'s. To do this we start with a lifting $\tilde x_0$
of the basepoint $x_0$ and we lift the other cells one after each other,
taking into account the relations in $\pi_1(E(K)$ and making sure that the union
of cells already lifted is always connected. When a cell $c$ is reached for the first
time, its lifting is chosen arbitrarily and declared to be $\tilde c^{(0)}$, but
its boundary will involve in general $\tilde d^{(n)}$'s with $n\neq 0$.
Once the lifted cellularization is known, it is a simple matter
to determine preferred cell liftings: since the support of $s'(P)$ consists
of 3 spiders, we only need to choose liftings of the 3 heads and then lift
the legs.

Carrying out the computations we have explicitly found the algebraic complex
with coefficients in $\mz[t^{\pm1}]$, and the preferred generators of the
4 moduli appearing. Then, using Maple, we have checked that indeed
the complex is acyclic, and we have computed its torsion as follows:
$$\tau^\varphi(E(K),\xi)=\pm t^{-1}.$$
Note that as an application of Proposition~\ref{wind:sensitive}, by adding curls, we can 
easily construct a family $\{K_n\}$ of pseudo-Legendrian knots such that
$\tau^\varphi(E(K_n),\xi_n)=\pm t^{n}$.

\section{Main proofs}\label{proofs}

In this section we provide all the proofs which we have
omitted in the rest of the paper. We will always refer to the statements
for the notation.

\dim{p:h:formula}
Let us first recall the classical Poincar\'e-Hopf formula, according to which if $v$ is a
vector field with isolated singularities on a manifold $M$, and $v$ points outwards on
$\partial M$ ({\em i.e.}~$\partial M$ is black), then the sum of the indices of all 
singularities is $\chi(M)$. Assume now that $v$ has isolated singularities and on 
$\partial M$ it is compatible with a pattern $\pp=(W,B,V,C)$. We claim that 
if $\cc$ is a cellularization of $M$ suited to $\pp$ we have:
\begin{equation}\label{ph:formula}
\sum_{x\in{\rm Sing}(v)}\index_x(v)=\chi(M)-\sum_{\sigma\in\cc,\ \sigma\subset W\cup V}
\index(\sigma).
\end{equation}
This formula is enough to prove the statement: if a non-singular field $v$ compatible
with $\pp$ exists then the left-hand side of~\ref{ph:formula} vanishes,
and the right-hand side of~\ref{ph:formula} equals the obstruction of the statement.
On the other hand, if the obstruction vanishes, then one can first consider a singular
field compatible with $\pp$, then group up the singularities in a ball, and remove them.

\begin{figure}
\centerline{\psfig{file=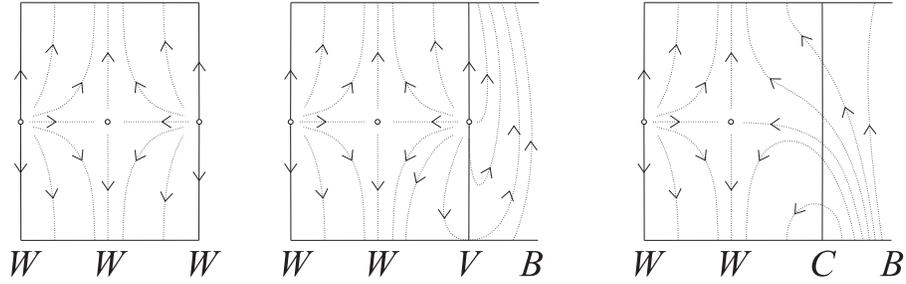,width=12cm}}
\caption{\label{two:obstr}Extension of a field to the collared manifold: dimension 2}
\end{figure}

To prove~\ref{ph:formula} we consider the manifold $M'$ obtained by attaching
a collar $\partial M\times[0,1]$ to $M$ along $\partial M=\partial M\times\{0\}$.
Of course $M'\cong M$. We will now extend $v$ to a field $v'$ on $M'$ with the property that
$v'$ points outwards on $\partial M'$, and in $\partial M\times(0,1)$ the field
$v'$ has exactly one singularity for each cell $\sigma\subset W\cup V$, with index
$\index(\sigma)$. An application of the classical Poincar\'e-Hopf formula then
implies the conclusion. The construction of $v'$ is done cell by cell. We first
show how the construction goes in dimension 2, see Fig.~\ref{two:obstr}.

\begin{figure}
\centerline{\psfig{file=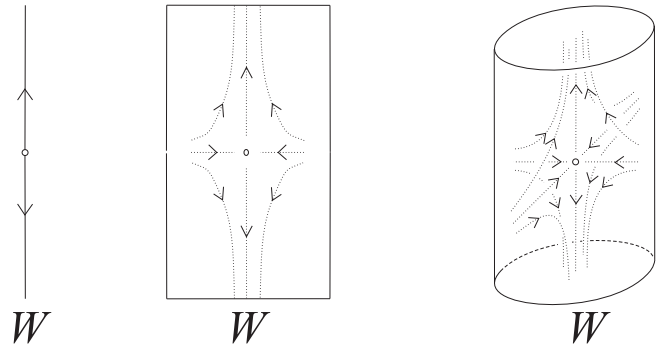,width=9cm}}
\caption{\label{three:obstr:1}Extension of a field to the collared manifold: white
cells in dimension 3}
\end{figure}

\begin{figure}
\centerline{\psfig{file=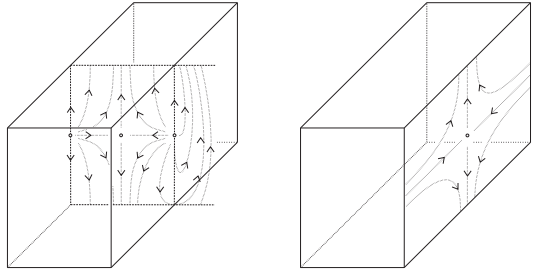,width=7cm}}
\caption{\label{three:obstr:2}Extension of a field to the collared manifold: convex edge in
dimension 3}
\end{figure}

For the 3-dimensional case, it is actually convenient to choose a cellularization $\cc$
of special type. Namely, we assume that $\cc\ristr{\partial M}$ 
consists of rectangles and triangles, each rectangle having exactly one edge on $V\cup C$,
and the union of rectangles covering a neighbourhood of $V\cup C$. We suggest 
in Fig.~\ref{three:obstr:1} how to 
define $v'$ on $\sigma\times[0,1]$ for $\sigma\subset W$ of dimension 0, 1 and 2 respectively.
By the choices we have made the situation near $\partial W$ contains the 2-dimensional
situation as a transversal cross-section, and it is not too difficult to extend $v'$ further
and check that indices of the singularities are as required. As an example, we suggest in
Fig.~\ref{three:obstr:2}
how to do this near a convex edge.\finedim{p:h:formula}

\begin{figure}
\centerline{\psfig{file=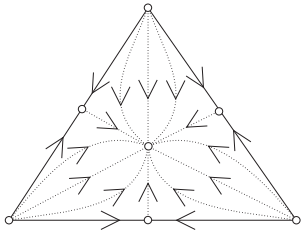,width=4cm}}
\caption{\label{elem:field}The singular field $w_\ss$ on a 2-simplex}
\end{figure}

\dim{reconstruction:statement}
Our proof follows the scheme given by Turaev in~\cite{turaev:Euler}, with some 
technical simplifications and some extra difficulties due to the tangency circles.
We first recall that if $\ss$ is a (smooth) triangulation of a manifold $N$, a (singular)
vector field $w_\ss$ on $N$ can be defined by the requirements 
that: (1) each simplex is a union of orbits; (2)
the singularities are exactly the barycentres of the
simplices; (3) barycentres of higher dimensional simplices are
more attractive that those of lower dimensional simplices.
More precisely, each orbit (asymptotically) goes from a barycentre $p_\sigma$
to a barycentre $p_{\sigma'}$, where $\sigma\subset\sigma'$.
It is automatic that $\index_{p_\sigma}(w_\ss)=\index(\sigma)$.
See Fig.~\ref{elem:field} for a description of $w_\ss$ on a 2-simplex of $\ss$.

Let us consider now a triangulation $\calt$ of $M$, and let us choose a 
representative $z$ of the given $\xi\in\eulc(M,\pp)$ as in
Proposition~\ref{combin:prop:statement}(\ref{bary:sub:point}).
We consider now the manifold $M$ obtained by attaching $\partial M\times[0,\infty)$
to $M$ along $\partial M=\partial M\times\{0\}$. Note that $M'\cong\interior(M)$.
Moreover $\calt$ extends to a ``triangulation'' $\calt'$ of $M'$, where on $M\times(0,\infty)$
we have simplices with exactly one ideal vertex, obtained by taking cones
over the simplices in $\partial M$ and then removing the vertex.
Even if $\calt'$ is not strictly speaking a triangulation, the construction of
$w_{\calt'}$ makes sense, because the missing vertex at infinity would be
a repulsive singularity anyway. We arrange things in such a way that if 
$\sigma\subset\partial M$ then the singularity in $\sigma\times(0,\infty)$
is at height $1$, so it is $p_\sigma\times\{1\}$.

We will define now a smooth function $h:\partial M\to (0,\infty)$ and set $M_h=M\cup\{(x,t)\in\partial M\times[0,\infty):\ t\leq h(x)\}$, in such a way that
$w_{\calt'}$ is non-singular on $\partial M_h$, and, modulo the natural homeomorphism
$M\cong M_h$, it induces on $\partial M_h$ the desired boundary pattern $\pp$.
Later we will show how to use $z$ to remove the singularities of $w_{\calt'}$ on $M_h$.

\begin{figure}
\centerline{\psfig{file=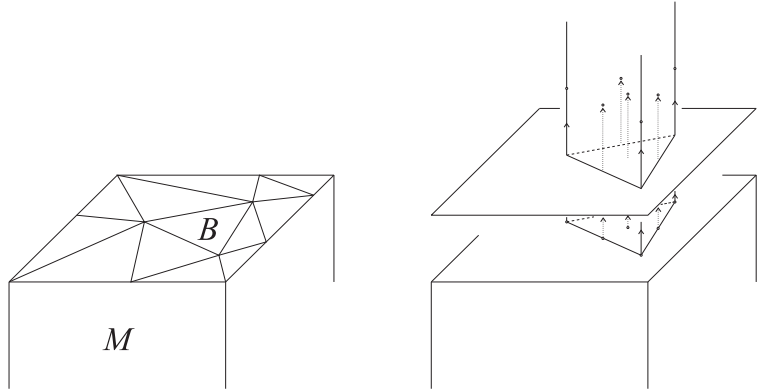,width=10cm}}
\caption{\label{blacktria}Where $h=1/2$ the field points outwards }
\end{figure}

\begin{figure}
\centerline{\psfig{file=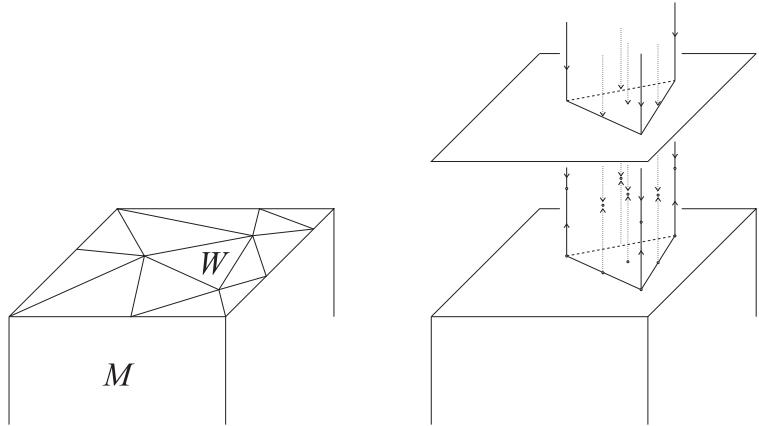,width=10cm}}
\caption{\label{whitetria}Where $h=2$ the field points inwards}
\end{figure}

\begin{figure}
\centerline{\psfig{file=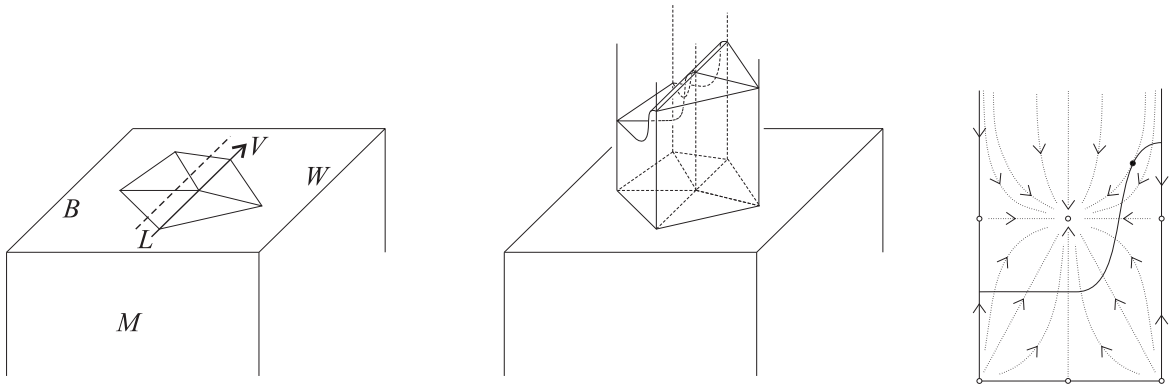,width=15.5cm}}
\caption{\label{convline}On $V$ the field has convex tangency}
\end{figure}

\begin{figure}
\centerline{\psfig{file=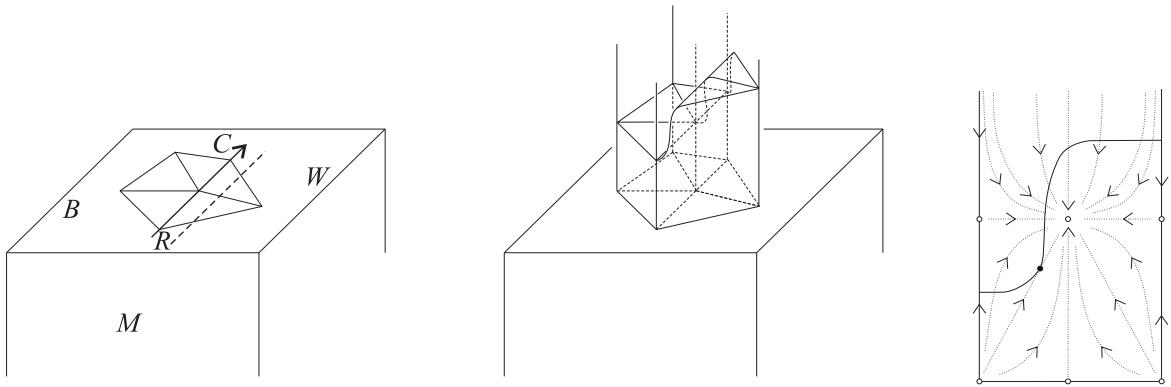,width=15.5cm}}
\caption{\label{concline}On $C$ the field has concave tangency}
\end{figure}

To define the function $h$ we consider a (very thin) left half-collar $L$ of $V$ on 
$\partial M$ and a right half-collar $R$ of $C$. Here ``left'' and ``right'' refer
to the natural orientations of $\partial M$ and of $V\cup C$.
Note that $L\subset B$ and $R\subset W$.
Now we set $h\ristr{B\setminus L}\equiv 1/2$, and 
$h\ristr{W\setminus R}\equiv 2$. Figures~\ref{blacktria} and~\ref{whitetria}
respectively show that away from $V\cup C$ indeed the pattern of $w_{\calt'}$ on
$\partial M_h$ is as required. Now we identify $L$ to $V\times[-1,0]$ and $R$ to
$C\times[0,1]$, and we define $h(x,s)=f(s)$ for $(x,s)\in V\times[-1,0]$
and $h(x,s)=f(s-1)$ for $(x,s)\in C\times[0,1]$, where $f:[-1,0]\to[1/2,2]$ is
a smooth monotonic function with all the derivatives vanishing at $-1$ and $0$.
Instead of describing $f$ explicitly we picture it and show that also near
$V\cup C$ the pattern is as required. This is done near $V$ and $C$ respectively 
in Figg.~\ref{convline}
and~\ref{concline}.
In both pictures we have only considered a special configuration for the 
triangulation on $\partial M$, and we have refrained from picturing the orbits
of the field in the 3-dimensional figure. Instead, we have separately shown
the orbits on the vertical simplices on which the value of $h$ changes.

The conclusion is now exactly as in Turaev's argument, so we only give a sketch.
The chosen representative $z$ of $\xi\in\eulc(M,\pp)$ can be described as an integer linear
combination of orbits of $w_{\calt'}$, which we can describe as segments
$[p_\sigma,p_{\sigma'}]$ with $\sigma\subset\sigma'$. Now we consider the chain
\begin{equation}\label{complete:chain:to:collar}
z'=z-\sum_{\sigma\subset W\cup V}\index(\sigma)\cdot p_\sigma\times[0,1].
\end{equation}
By definition of $h$ we have that $z'$ is a 1-chain in $M_h$, and 
$\partial z'$ consists exactly of the singularities of $w_{\calt'}$ contained in $M_h$, 
each with its index. For each segment $s$ which appear in $z'$
we first modify $w_{\calt'}$ to a field which is ``constant''
on a tube $T$ around $s$, and then we modify the field again
within $T$, in a way which depends on the coefficient of $s$ in $z'$.
The resulting field has the same singularities as $w_{\calt'}$, but 
one checks that these singularities can be removed by 
a further modification supported within small balls centred at the singular points.
We define $\Psi(\xi)$ to be the class in $\euls(M,\pp)$ of this final field.
Turaev's proof that $\Psi$ is indeed well-defined and $H_1(M;\mz)$-equivariant
applies without essential modifications.\finedim{reconstruction:statement}

\begin{rem}\label{good:cell:ok}
{\em In the previous proof we have defined $\Psi$ using triangulations, in order to
apply directly Turaev's technical results (in particular, invariance under subdivision).
However the geometric construction makes sense also for cellularizations $\cc$ more
general than triangulations, the key point being the possibility of defining a field
$w_{\cc}$ satisfying the same properties as the field defined for triangulations.
This is certainly true, for instance, for cellularizations $\cc$ of $M$ induced by
realizations of $M$ by face-pairings on a finite number of polyhedra, assuming
that the projection of each polyhedron to $M$ is smooth.}
\end{rem}

\dim{diagram:commutes}
To help the reader follow the details, we first outline the scheme of the proof:
\begin{enumerate}
\item By identifying $M$ to a collared copy of itself, we choose a representative $z$
of the given $\xi\in\eulc(M,\pp)$ such that the extra terms added to define $\Thetac(\xi)$
cancel with terms already appearing in $z$. (We know {\em a priori} that this happens
at the level of boundaries, but it may well not happen at the level of $1$-chains.)
\item We apply Remark~\ref{good:cell:ok} and choose a cellularization of $M$ in which
it is particularly easy to analyze $\Psi(\xi)$ and $\Psi(\Thetac(\xi))$, both constructed
using the representative $z$ already obtained.
\end{enumerate}

\noindent We consider a cellularization $\cc$ of
$M$ satisfying the same assumptions on $\partial M$ as those considered in the proof of
Proposition~\ref{p:h:formula}, namely $C\cup V$ is surrounded on both sides
by a row of rectangular tiles, and the other tiles are triangular. 
We denote by $\gamma_1,\dots,\gamma_n$ the segments in $C$, oriented as $C$.

Let us consider a representative $z$ relative to $\cc$ of the given $\xi\in\eulc(M,\pp)$.
We construct a new copy $M_1$ of $M$ by attaching $\partial M\times[-1,0]$ to
$M$ along $\partial M=\partial M\times\{-1\}$, and we extend $\cc$ to $\cc_1$ by 
taking the product cellularization on $\partial M\times[-1,0]$. We define a new chain as
\begin{eqnarray*}
z_1 & = & z+\sum_{\sigma\subset B}\index(\sigma)\cdot p_\sigma\times[-1/2,0]
-\sum_{\sigma\subset W\cup V}\index(\sigma)\cdot p_\sigma\times[-1,-1/2]\\
& & + \sum_{j=1}^n\Big(\gamma_j\ristr{[1/2,1]}\times\{-1/2\}-
\gamma_j\ristr{[1/2,1]}\times\{0\}\Big).
\end{eqnarray*}
Note that $z_1$ is an Euler chain in $M_1$ with respect to $\cc_1$. Consider the natural
homeomorphism $f:M\to M_1$ and the class
$$a=\alphac(f_*(\xi),[z_1])\in H_1(M_1;\mz)$$
which may or not be zero. Since the inclusion of $M$ into $M_1$ is an isomorphism
at the $H_1$-level, $a$ can be represented by a $1$-chain in $M$, so $z_1$ can
be replaced by a new Euler chain $z_2$ such that $[z_2]=f_*(\xi)$ and $z_2$
differs from $z_1$ only on $M$.

Renaming $M_1$ by $M$ and $z_2$ by $z$ we have found a representative $z$ of $\xi$
such that $z=z_\theta+\sum_{j=1}^n\gamma_j\ristr{[1/2,1]}$, where $z_\theta$ is a sum
of simplices contained in $B\cup\interior M$. Note that of course $\Thetac(\xi)=[z_\theta]$.
To conclude the proof we need now to analyze $\Psi(\xi)$, constructed using $z$, and
$\Psi(\Thetac(\xi))$, constructed using $[z_\theta]$, and show that
$\Psi(\Thetac(\xi))=\Thetas(\Psi(\xi))$. By construction $\Psi(\xi)$ and 
$\Psi(\Thetac(\xi))$ will only differ near $C$, and we concentrate on one component
of $C$ to show that the difference is exactly (up to homotopy) as in the definition of
$\Thetas$, {\em i.e.}~as in Fig.~\ref{conc:to:conv}.

The difference between $\Psi(\xi)$ and $\Psi(\Thetac(\xi))$ is best visualized on
a cross-section of the form $C\times[0,\infty)$. We leave to the reader to analyze
the complete 3-dimensional pictures. To understand the cross-section, we 
follow the various steps in the proof of Theorem~\ref{reconstruction:statement}.

\begin{figure}
\centerline{\psfig{file=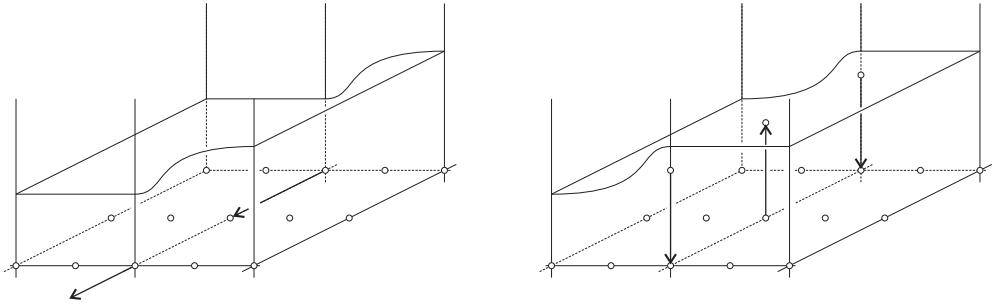,height=4cm}}
\caption{\label{comm:diag:1}Local difference between $z'$ (left) and $z'_\theta$ (right)}
\end{figure}

The first step in the definition of $\Psi(\xi)$ (respectively, $\Psi(\Thetac(\xi))$)
consists in choosing the height function $h$ (respectively, $h_\theta$)
and replacing the chains $z$ (respectively, $z_\theta$)
by a chain $z'$ (respectively, $z'_\theta$)
as in formula~(\ref{complete:chain:to:collar}).
This is done in Fig.~\ref{comm:diag:1}
where only the difference between the chains is shown.

\begin{figure}
\centerline{\psfig{file=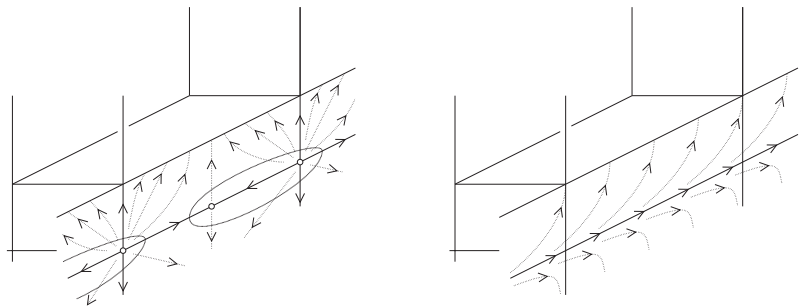,height=4cm}}
\caption{\label{comm:diag:2}Construction of $\Psi(\xi)$ on $C\times[0,\infty)$.
On the left we show $w_\cc$ and the zones where it must be modified.}
\end{figure}

\begin{figure}
\centerline{\psfig{file=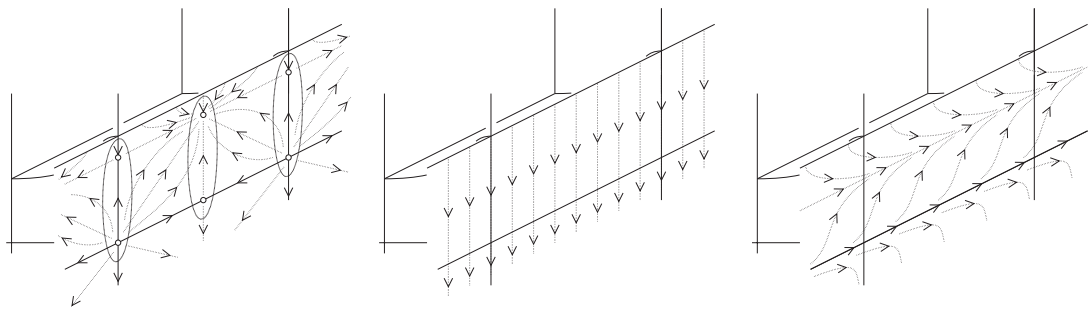,height=4cm}}
\caption{\label{comm:diag:3}Construction of $\Psi(\Thetac(\xi))$ on $C\times[0,\infty)$}
\end{figure}

To conclude we must modify the field $w_\cc$ within a small neighbourhood of 
the support of $z'$ and $z'_\theta$. This is done in Figg.~\ref{comm:diag:2}
and~\ref{comm:diag:3} respectively. The rightmost picture in Fig.~\ref{comm:diag:3} 
is obtained by homotopy on the previous one. The representatives of $\Psi(\xi)$ and
$\Psi(\Thetac(\xi))$ can be compared directly, and indeed they differ by a curve
parallel to $C$ and directed consistently with $C$, so 
$\Psi(\Thetac(\xi))=\Thetas(\Psi(\xi))$.\finedim{diagram:commutes}

We give now the proof omitted in Section~\ref{link:section}.

\dim{wind:sensitive}
Let us first note that the comparison class which we must show to be $[m]$ is independent
of $f$ by Proposition~\ref{collar:trivial}. We
will give two completely independent (but somewhat sketchy) proofs that this
class is indeed $[m]$.

\begin{figure}
\centerline{\psfig{file=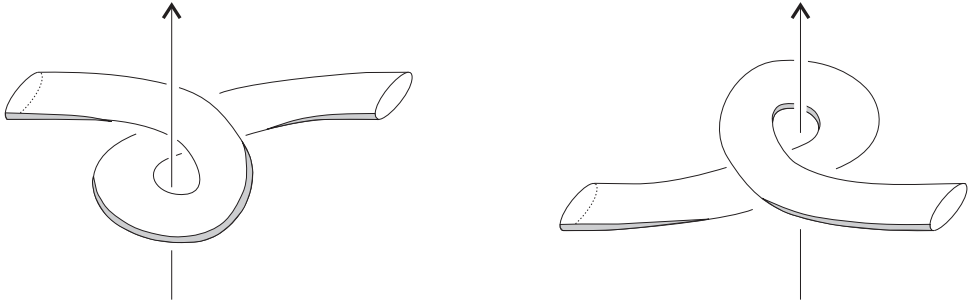,height=4cm}}
\caption{\label{curls}Differently curled tubes in the vertical field.}
\end{figure}

\begin{figure}
\centerline{\psfig{file=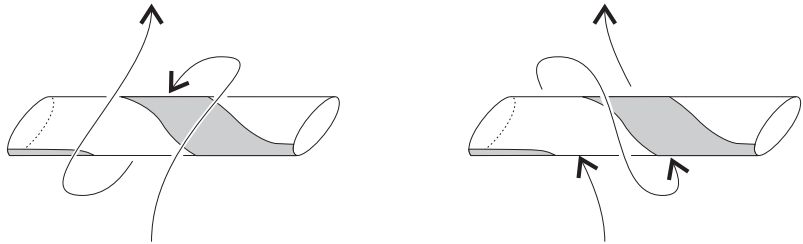,height=3.2cm}}
\caption{\label{nocurls}Straightened curls.}
\end{figure}

For a first proof, instead of comparing a ``straight'' knot with one with two curls, 
we compare two knots with one curl, chosen so that the framing is the same but the 
winding number is different. This is of course equivalent. The two knots are shown in
Fig.~\ref{curls} as thick tubes, together with one specific orbit of the
field they are immersed in. The resulting bicoloration on the boundary
of the tubes is also outlined.
To compare the curls we isotope the bicolorated tubes to the same
straight tube, and we show how the orbit of the field is transformed under
this isotopy. This is done in Fig.~\ref{nocurls}.
Also from this very partial picture it is quite evident that the resulting fields
wind in opposite directions around the tube. A more accurate picture
would show that the difference is actually a meridian of the tube.

Another (indirect) proof goes as follows. Note first that the comparison class which we
must compute certainly is a multiple of $[m]$, say $k\cdot[m]$. Note also that
$k$ is independent of the ambient manifold $(M,v)$. Moreover, by symmetry, we
will have $\alpha(\xi(v,K_0),f_*(\xi(v,K_{-1})))=-k\cdot[m]$ if $K_{-1}$ is
obtained by locally adding a double curl with opposite winding number.

We take now $M$ to be $S^3$, with the field $v$ carried by the abalone $P$
as in Section~\ref{exa:section}, and $K$ to be a trivial knot contained in the 
``smaller'' disc of the $P$. We apply Proposition~\ref{link:family:prop}
to find another pseudo-Legendrian knot $K'$ in $(S^3,v)$ such that 
$\alpha(\xi(v,K),\xi(v,K'))=[m]$, where by simplicity we are omitting the
framed-isotopies necessary to compare these classes. As already remarked in
Section~\ref{spines:section}, we can assume that $K'$ has a $\cont^1$-projection
on $P$. If one examines $P$ carefully one easily sees that $K'$ can actually be slid
over $P$ to lie again in the small disc of $P$. Now $K'$ is a planar projection of the
trivial knot, so through Reidemeister moves of types II and III,
which correspond to isotopies throgh knots transversal to $v$, it can be transformed
into a projection which differs from the trivial one only for a finite (even) number
of curls. Summing up, we have a knot $K'$ such that $\alpha(\xi(v,K),\xi(v,K'))=[m]$
and $K'$ differs from $K$ only for a finite number of transformations of the form
$K\mapsto K_1$ or $K\mapsto K_{-1}$. This shows that $[m]$ is a multiple of 
$k\cdot[m]$, so $k=\pm1$.\finedim{wind:sensitive}

We conclude the paper by establishing
the only statement given in Section~\ref{spines:section} and
not proved there. As above, we do not recall all the notation.

\dim{spider:structure:teo}
We fix $P$ and set $s''=s''(P)$, $\hatv=\hatv(P)$.
Using Remark~\ref{good:cell:ok} we see that the construction of
$\Psi([s''])$ explained in the proof of Theorem~\ref{reconstruction:statement}
can be directly applied to the cellularization $\hattt=\hattt(P)$ of $\hatM$.
Recall that this construction requires identifying $\hatM$ to a collared copy
of itself, and extending $s''$ to a chain $s'''$ whose boundary consists precisely
of the singularities of a field $w$. A representative of $\Psi([s''])$ is then
obtained by applying to $w$ a certain desingularization procedure.
This desingularization is supported in a
neighborhood of $s'''$, and one can easily check that
the connected components of the support of $s'''$ (denoted henceforth by $S$)
are actually contractible.
Therefore, {\em any} desingularization of $w$ supported in a
neighbourhood of $s'''$ will give a representative of $\hatv$.
We will prove the desired formula $\Psi([s''])=[\hatv]$
by exhibiting one such desingularization which is nowhere antipodal to $\hatv$.
In our argument we will always neglect the contraction of $\stwotriv$ which maps
$M$ onto $\hatM$. (The desired formula actually holds at the level of $M$, and it
easily implies the formula for $\hatM$.) 

\begin{figure}
\centerline{\psfig{file=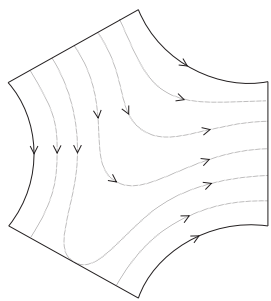,width=3.5cm}}
\caption{\label{hexagon:1}The field $\hatv$ on a hexagon}
\end{figure}

\begin{figure}
\centerline{\psfig{file=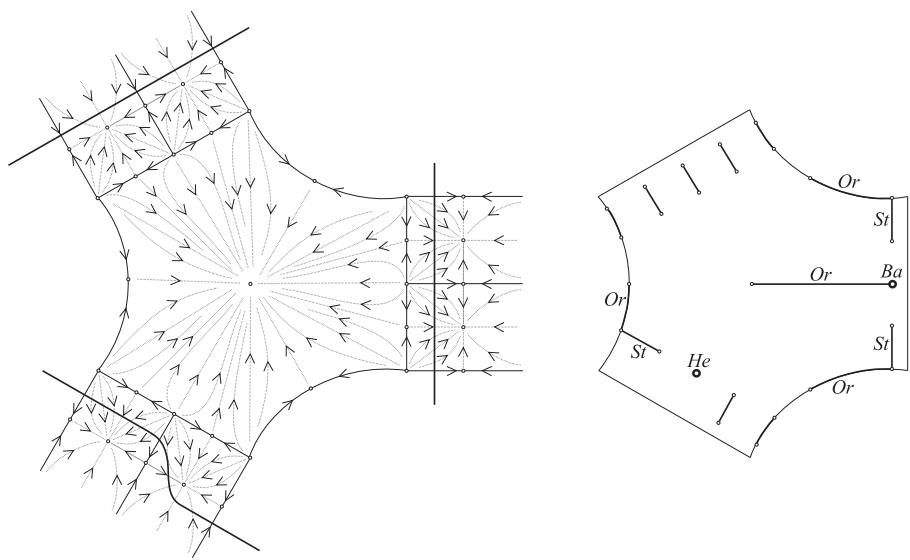,width=12cm}}
\caption{\label{hexagon:2}The field $w$ and the trace of $S$ on a hexagon}
\end{figure}

By the above observations, the following claims easily imply the conclusion of the proof:
\begin{enumerate}
\item\label{claim:antip} The set of points where $w$ is antipodal to 
$\hatv$ is contained in $S$.
\item\label{claim:desing} If $S_0$ is a component of
$S$ then $w$ can be desingularized within a neighbourhood of
$S_0$ to a field which is not antipodal to $\hatv$ in the neighbourhood.
\end{enumerate}
\noindent We prove claim~\ref{claim:antip} by first noting that the cells dual
to those of $P$ are unions of orbits of both $w$ and $\hatv$. So we can analyze cells
separately. We do this explicitly only for 2-dimensional cells, leaving to 
the reader the other cases. In Fig.~\ref{hexagon:1}
we describe $\hatv$. In the left-hand side of Fig.~\ref{hexagon:2}
we describe $w$ on the collared hexagon. In the right-hand side of the same figure
we only show the singularities of $w$ on the renormalized hexagon,
and the intersection of $S$ with the hexagon. In this figure the 7 short segments come
from $s'''-s''$; the other bits of $S$ have been labeled by `Or', `St', `Ba' or `He'
to indicate that they come from
orbits of $\hatv$, stars, bi-arrows or half-edges.

\begin{figure}
\centerline{\psfig{file=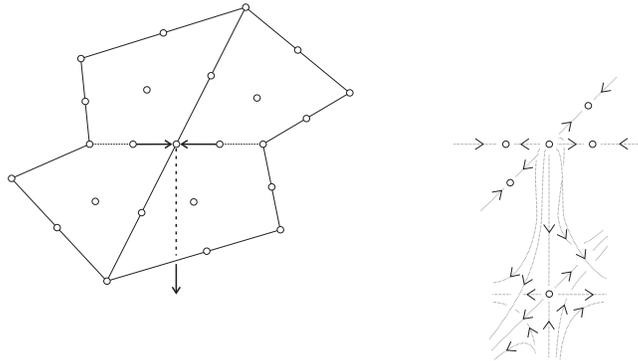,width=8.5cm}}
\caption{\label{fix:biarrow}An enhanced bi-arrow and the field $w$ near it}
\end{figure}

This proves claim~\ref{claim:antip}. Comparing Fig.~\ref{hexagon:2}
with Fig.~\ref{hexagon:1}, and carrying out the same analysis for
3-cells, one actually shows also claim~\ref{claim:desing} for 
components $S_0$ coming from $s'''-s''$. Components of $S$ other than these can
be described in one of the following ways:
\begin{itemize}
\item[(a)] an orbit of $\hatv$ emanating from a vertex of $P$;
\item[(b)] a half-edge of $C$;
\item[(c)] a bi-arrow together with an orbit of $\hatv$ emanating from
the centre of an edge of $P$ and reaching the centre of the bi-arrow;
\item[(d)] a star together with an orbit of $\hatv$ emanating from
the centre of a disc of $P$ and reaching the centre of the star.
\end{itemize}
\noindent All cases can be treated with the same method, we only do case (c). 
Figure~\ref{fix:biarrow}
shows the component placed so that $\hatv$ can be thought of as the 
constant vertical field pointing upwards, and the field $w$ near the component.
The conclusion easily follows.
\finedim{spider:structure:teo}

\vspace{1cm}

\hspace{8cm} benedett@dm.unipi.it

\hspace{8cm} petronio@dm.unipi.it

\hspace{8cm} Dipartimento di Matematica

\hspace{8cm} Via F.~Buonarroti, 2

\hspace{8cm} I-56127, PISA (Italy)

\end{document}